\documentclass[11 pt]{amsart}
\usepackage{amssymb}
\usepackage{amsmath}
\usepackage{amsfonts}
\usepackage{graphicx}
\usepackage{amsthm}
\usepackage{enumerate}
\usepackage[mathscr]{eucal}
\usepackage{verbatim}
 \theoremstyle{plain}
\newtheorem{theorem}{Theorem}[section]
\newtheorem{lemma}[theorem]{Lemma}
\newtheorem{proposition}[theorem]{Proposition}

\numberwithin{equation}{section}
\theoremstyle{plain}

\numberwithin{equation}{section}
\theoremstyle{remark}

\def\conv{\text{conv}}

\def\tcr{{t_{\text{cr}}}}

\def\intslash{\rlap{\kern  .32em $\mspace {.5mu}\backslash$ }\int}
\def\qsl{{\rlap{\kern  .32em $\mspace {.5mu}\backslash$ }\int_{Q_x}}}

\def\vpi{\varpi}
\def\vth{\vartheta}

\def\emph#1{{\it #1 }}

\def\cf{{\it cf}}

\def\supp{{\text{\rm supp }}}

\def\inn#1#2{\langle#1,#2\rangle}

\def\noi{\noindent}

\def\meas{{\text{\rm meas}}}

\def\lc{\lesssim}

\def\eps{\varepsilon}

\def\ka{\kappa}

\def\la{\lambda}
\def\La{\Lambda}

\def\Om{\Omega}
\def\fA{{\mathfrak {A}}}

\def\fM{{\mathfrak {M}}}

\def\fS{{\mathfrak {S}}}

\def\fa{{\mathfrak {a}}}

\def\bbN{{\mathbb {N}}}

\def\bbR{{\mathbb {R}}}

\def\bbZ{{\mathbb {Z}}}
\def\cA{{\mathcal {A}}}
\def\cB{{\mathcal {B}}}

\def\cE{{\mathcal {E}}}

\def\cK{{\mathcal {K}}}

\def\cM{{\mathcal {M}}}

\def\cR{{\mathcal {R}}}

\def\cT{{\mathcal {T}}}

\def\cV{{\mathcal {V}}}

\def\cZ{{\mathcal {Z}}}

\def\tx{\tilde x}
\def\ts{\tilde s}

\def\tw{\tilde w}

\font \roman = cmr10 at 10 true pt

\def\be#1{\begin{equation}\label{#1}}
\def\ee{\end{equation}}
\def\bas{\begin{align*}}
\def\eas{\end{align*}}
\def\bi{\begin{itemize}}
\def\ei{\end{itemize}}

\def\supp{{\hbox{\roman supp}}}

\def\eps{\varepsilon}

\def\emph#1{{\it #1}}
\def\textbf#1{{\bf #1}}
\def\intslash{\rlap{\kern  .32em $\mspace {.5mu}\backslash$ }\int}
\def\qsl{{\rlap{\kern  .32em $\mspace {.5mu}\backslash$ }\int_{Q_x}}}

\begin{document}
\title
[Fourier restriction for space curves]
{Restriction of Fourier transforms to curves and related oscillatory integrals}

\author[]
{Jong-Guk Bak \ \ Daniel M.  Oberlin \  \  Andreas Seeger}

\address {J. Bak \\Department of Mathematics and the Pohang
  Mathematics Institute\\ Pohang University of Science and Technology
\\
Pohang 790-784, Korea}
\email{bak@postech.ac.kr}

\address
{D. M.  Oberlin \\
Department of Mathematics \\ Florida State University \\
 Tallahassee, FL 32306}
\email{oberlin@math.fsu.edu}

\address{A. Seeger   \\
Department of Mathematics\\ University of Wisconsin-Madison\\Madison,
WI 53706, USA}
\email{seeger@math.wisc.edu}


\thanks{J.B. was supported in part by  grant
R01-2004-000-10055-0 of the Korea Science and Engineering Foundation,
D.O.  was supported in part by NSF grant DMS-0552041, and
A.S. was supported in part by NSF grant DMS-0200186.}

\begin{abstract}
We prove  sharp
endpoint results for the Fourier restriction operator
associated to nondegenerate curves in $\Bbb R^d$, $d\ge 3$,
and related  estimates for
oscillatory integral operators. Moreover, for some larger classes of curves
in $\Bbb R^d$ we obtain  sharp  uniform $L^p\to L^q$
bounds with respect to affine arclength measure, thereby resolving a problem of Drury and Marshall.

\end{abstract}

\maketitle

\section{Introduction}

\label{intro}

For a Schwartz function $f$ defined on $\bbR^d$, $d\ge 2$
 consider the restriction of its Fourier transform
to the curve $t\mapsto (t,t^2,\dots, t^d)$
{\it i.e.},
\[
\cR f(t)= \widehat f(t,t^2,\dots, t^d).
\]
The problem of $L^p(\bbR^d)\to L^q(\bbR)$ bounds is
 understood, namely $\cR$ is bounded from $L^p(\bbR^d)$ to $L^q(\bbR)$
if and only if
\begin{equation}\label{pqtypeset}
1\le p<p_d:= \frac{d^2+d+2}{d^2+d}
 \quad\text{ and } p'=\frac{d(d+1)}{2}  q;
\end{equation}
likewise,  if $I$ is a compact interval, then,
for the same range of $p$,
$\cR$ is bounded from $L^p(\bbR^d)$ to $L^q(I)$
when
$p'\ge d(d+1)q/2$.

In two dimensions the sharp boundedness  result is due to Zygmund \cite{Z}
who extended earlier
work by Fefferman and Stein (\cite{F}) (see also \cite{CS}, \cite{H}
for estimates on more general oscillatory integral operators).
Initial results in higher dimensions
for the smaller range
$1\le p<(d^2+2d)/(d^2+2d-2)$ are  due to Prestini \cite{Pr},
 with strict inequality
$p'>
 d(d+1)q/2$
 for the local result. For the same range
 of  $p$, Christ \cite{Ch} showed boundedness on the edge
$p'=
 d(d+1)q/2$.
The full range
\eqref{pqtypeset} has been obtained by Drury \cite{D}.
The necessity of the condition $p'\ge d(d+1)q/2$  follows by scaling and
the necessity of the
condition $1\le p<p_d$ follows from work by Arkhipov, Chubarikov
and Karatsuba \cite{aktdokl} ({\it}\cf. also  \cite{mock},
\cite{akt} and \cite{bggist}).

The first problem considered in this paper is
what can be said about estimates for
 the endpoint $p=q=p_d$.
In two dimensions $p_2=4/3$ and
Beckner,  Carbery,  Semmes and   Soria \cite{bcss} showed by a Kakeya set argument that the restricted weak type $4/3$ estimate  fails; in fact
$\cR$ does not even map $L^{4/3,1}(\bbR^2)$ to  $L^{1,\infty}(I)$.
 Using a result by Keich \cite{Ke} this can be further quantified
for functions
supported  in $D_R:=\{x\in\bbR^2:|x|\le R\}$, namely one gets for  large $R$
the lower bound in the equivalence
\begin{equation} \label{2dinDR}
\sup \big\{ \|\cR f\|_{L^{1,\infty}(I)}: \|f\|_{L^{4/3,1}(\bbR^2)}\le 1, \,
\supp(f)\subset D_R \big \} \approx  [\log R]^{1/4}.
\end{equation}
The upper bounds can be deduced from
the method or results in  \cite{H}
 and indeed the
analogue of \eqref{2dinDR}
for the  $L^{4/3}(\bbR^2)\to L^{4/3}(I)$ inequality
 holds as well, for functions supported in $D_R$ (see also  \cite{tomas}
for a related  estimate).

In higher dimensions the arguments  by Drury are
not promising to obtain similar endpoint bounds. He obtained
 his result by an elegant  iteration  procedure where it is shown
that if the   $L^{a}\to L^{b}$ inequality  for the adjoint operator
holds on the  critical edge  for  some
range $b>b_0>q_d:=p_d'=(d^2+d+2)/2$, then it also holds for $b>\rho(b_0)$
where $b_0>\rho(b_0)>q_d$ and the sequence recursively defined by
$b_{i+1}=\rho(b_i)$ is decreasing and converges to $q_d$.
We note that the constants in the estimates
increase exponentially with the number of iterations, so that
a sensible endpoint bound for $\cR$ on  functions in $L^{p_d,1}$ seems out of
reach with this method.
Given also the lower bounds
 \eqref{2dinDR}
 in two dimensions, it is  somewhat
surprising that the restricted weak type endpoint bound does hold
 in three and higher dimensions; in fact the   better  {\it restricted strong type} estimate is
true:

\begin{theorem}\label{restr}
Let $d\ge 3$ and  $p_d=\tfrac{d^2+d+2}{d^2+d}$.
Then
\[ \cR: L^{p_d,1}(\bbR^d) \to L^{p_d}(\bbR)\]
is bounded.
\end{theorem}
Note that all $L^p\to  L^q$ inequalities for $\cR$ can be deduced by
interpolation  with trivial
estimates for $L^1$ functions.
Arguments in \cite{akt} or \cite{bggist} show that the
source space
$L^{p_d,1}$ cannot be replaced by a larger Lorentz space.
The argument in \S\ref{sharpness} below also shows that
 the target space $L^{p_d}$ cannot be replaced by a smaller Lorentz
space. It would be interesting  to investigate  whether the validity of the
endpoint bounds in higher dimensions
has implications to some integral geometric or Kakeya type
problems.

The estimate for $\cR$ is deduced
from an estimate for the adjoint operator which after some rescaling leads to a problem about more general oscillatory integral operators
given by
\begin{equation}\label{Tla}
T_\la f(x)= \int e^{i\la \phi(x,t)} a(x,t) f(t) dt
\end{equation}
where $\la \gg 1 $,
the amplitude $a$ is $C^\infty$ and  compactly supported in
$\Omega\times I \subset \bbR^d\times \bbR$
and  $\phi$ is a real valued phase function in $C^\infty (\Om\times I)$.
Following \cite{baklee} we impose the  curvature condition
that for each $x\in \Om$  the curve
$s\mapsto \nabla_x\phi(x,s)$ is nondegenerate, in the sense that
\begin{equation}
\det\big(  \partial_{t}(\nabla_{x}\phi)  ,\,\partial_{t}^{2}
(\nabla_{x}\phi)  ,\,\ldots,\,\partial_{t}^{n}(\nabla_{x}\phi)  \big)
 \neq0 \label{nondeg}
\end{equation}
in $\Omega\times I$.

\begin{theorem}\label{osc}
Suppose that $d\ge 3$ , $\la>0$,
and that
\eqref{nondeg} holds in $\Omega\times I$.
Let  $ q_d=\frac{d^2+d+2}{2}$.
Then
 \begin{equation}
\label{rleq}
\big\|T_\la \big\|_{L^{q_d}(I)\to  L^{q_d, \infty}(\bbR^d)}
\le C (1+\la)^{-d/q_d}.
\end{equation}
Moreover,
\begin{equation}\label{Lqbd}
\big\|T_\la \big\|_{L^{q_d}(I)\to  L^{q_d}(\bbR^d)} \le C \big(\log (2+\la)\big)^{1/q_d}
 (1+\la)^{-d/q_d}.
\end{equation}
\end{theorem}
Again
\eqref{rleq} and even the weaker restricted weak type inequality
 fail in two dimensions,
 by a Kakeya set argument (\cite{bcss}, \cite{Tao-CS}).
Calculations  with
 $f\in C^\infty_0$  and the phase for the extension operator,
$\phi(x,t)=-\inn{\gamma(t)}{x}$,
show
that \eqref{Lqbd} is sharp; moreover the space
$L^{q_d, \infty}(\bbR^d) $ in \eqref{rleq}
 cannot be replaced by a smaller Lorentz space
$L^{q_d,s}$, see  \cite{akt}, \cite{bggist}.
Finally we shall show in  \S\ref{sharpness} that for the case of the
extension operator
 in \eqref{rleq} $L^{q_d}$
 cannot be replaced by any Lorentz space
$L^{q_d,\rho}$ with $\rho>q_d$.
We point out an important aspect  of the proof of Theorem \ref{osc},
namely the idea that in order to prove a weak type $(q,q)$ bound for
very large $q$ one estimates
a  multilinear expression with many symmetries
on the space $L^{r,\infty}$, for some  $r<1$, and takes advantage of the
{\it $r$-convexity} of this space
(\cf.
\S\ref{prelim}).

\medskip

We now turn to classes of curves for which the
nondegeneracy condition is  not satisfied.  It has long been known that
restriction theorems such as Theorem \ref{restr}  hold under a finite type condition,
if $p$ is taken from a sufficiently small interval $[1, 1+\epsilon)$ with $\epsilon$ depending on the ``type'';
for the known results of this type  see the papers
by Sj\"olin \cite{sj} and Sogge \cite{so} in two dimensions and by
 Christ \cite{Ch} and  Drury and Marshall
\cite{DM1}, \cite{DM2} in higher dimensions.
Another direction that has been pursued is to prove
a sharp universal
 restriction theorem in the {\it full range} $p\in [1,p_d)$, on the
critical edge $1/q= d(d+1)/(2p')$. Now the standard  measure needs to be
 replaced by the {\it affine arclength measure}
given by $w(t) dt$ where
\begin{align}\label{wdef}
w(t)&\equiv w_\gamma(t)= |\tau(t)|^{\frac{2}{d(d+1)}},
\\
\label{taudef}
\tau(t)&\equiv \tau_\gamma(t) = \det
\begin{pmatrix} \gamma'(t),& \gamma''(t),&  \dots,& \gamma^{(d)}(t)
\end{pmatrix}.
\end{align}
The objective is then to prove   the endpoint inequality
\begin{equation} \label{affrestr}
\Big(\int_I |\widehat f (\gamma(t))|^q w(t) dt\Big)^{1/q} \le C\|f\|_p,
\quad p'= \tfrac{d(d+1)}{2}q,
\end{equation}
for all $p<p_d$, see  {\it e.g.} \cite{D2} for a general discussion.

Note that  the arclength measure is invariant
 under  reparametrization. Moreover, an important feature
of the inequality
\eqref{affrestr} on the critical edge
 is its invariance under general linear change of
 variables.

In two dimensions, inequality \eqref{affrestr}
has been proved  by Sj\"olin \cite{sj}
for large classes of  convex curves, see also  Oberlin \cite{ob}.
Moreover Drury and Marshall
\cite{DM2}
proved a positive result for  large classes of finite type curves in
higher dimensions, in the partial
range $p<\tfrac{d^2+2d}{d^2+2d-2}$ ({\it i.e.} $ p<15/13$ in three dimensions),
and Drury \cite{D2} improved the result for some curves
 in three dimensions, obtaining a better partial result for  $p<36/31$.

It is conceivable that the
inequality \eqref{affrestr} is true for all curves, as to our best
knowledge no counterexamples are known.
However, in  three and higher dimensions,
very few positive results have been known for all $p<p_d$.
For example given the family of monomial curves with
nonvanishing curvature, {\it i.e.}  $(t, t^2, t^\beta)$,   the only three cases
 for which
\eqref{affrestr}
has been  known in the full range $p<p_3=7/6$
are (i) the trivial case $\beta=2$ where the weight $w$ vanishes
identically,
(ii) the nondegenerate case $\beta=3$, and, (iii),  the curious  exceptional
 case $\beta=9$, which follows via a change of variables
 from an  estimate for ``rough'' nondegenerate curves, due to the
first  two authors  \cite{BO}, \cf. also Theorem
\ref{BOext}
below.

We prove sharp and uniform $L^p\to L^q$
estimates for all monomial curves.

\begin{theorem} \label{powerthm} For given real numbers  $a_1,\dots,a_d $, $d\ge 2$,  consider the curve
$$t\mapsto \gamma(t)=(t^{a_1},t^{a_2},\dots, t^{a_d}), \quad
0<t<\infty,$$ let $w$ be as in
\eqref{wdef}
 and suppose $1<p<p_d
=\tfrac{d^2+d+2}{d^2+d}$, and $p'= \tfrac{d(d+1)}{2}q$.
Then there is $C(p,d)<\infty$ so that for all $f\in L^p(\bbR^d)$
\begin{equation} \label{affrestr2}
\Big(\int_0^\infty
 |\widehat f (\gamma(t))|^q w(t) dt\Big)^{1/q} \le
C(p,d)\|f\|_{L^p(\bbR^d)}.
\end{equation}
\end{theorem}
It should be emphasized  that the finite constant $C(p,d)$ does not
depend on the choice of  exponents $a_1,\dots, a_d$.
Some related results for  classes  of `simple'
 curves $(t,t^2,\dots, t^{d-1},\phi(t))$,
but with possibly flat $\phi$  will
be  treated in a subsequent paper \cite{bos2}.

\medskip

Finally, it is natural to ask whether a restricted strong type estimate
 with respect  to affine arclength measure
holds at the endpoint $p=p_d$, $d\ge 3$,
 for some class of
``degenerate curves''. This remains largely open, and we have a positive result only
for special cases. We formulate such a result for
certain curves in $\Bbb R^3$; note that the critical exponent is $p_3=7/6$.

\begin{theorem}\label{BOext}
Let  $\gamma (t) = (t, t^\alpha, t^{\beta})$, and
$w(t)
=|\det( \gamma'(t),
\gamma''(t), \gamma'''(t))|^{1/6}.$ Then there is a universal
 constant $C$ so that
the inequality
\begin{equation}\label{globalwt}
\Big(\int_0^\infty \big| \widehat{f}(\gamma (t))\big |^{7/6} w(t) dt\Big)^{6/7}
\le  C \big\| f\big\|_{L^{7/6, 1} (\bbR^3)}.
\end{equation}
holds in the following two cases:

(i)  $\alpha+\beta=5$, $\alpha\notin (2,3)$,

(ii) $\beta=5\alpha-1$, $\alpha\notin (1/3, 1/2)$.
\end{theorem}
The proof of this result is a combination of the method in \S
\ref{lorentzbd} with ideas in \cite{BO}.

\medskip


{\it Structure of the paper:} In \S\ref{prelim} we discuss
preliminaries on Lorentz spaces, multilinear interpolation and the
Drury-Marshall bound on a class of multilinear operators
involving Vandermonde determinants (the proof is given in an
appendix). The weak type estimates for nondegenerate curves
are proved in \S\ref{lorentzbd}, and the strong type bound \eqref{Lqbd}
is proved in \S\ref{strongtypebound}. \S\ref{sharpness} contains a lower bound
for the norms of the  extension operator proving the sharpness of the
weak type $q_d$ bound.
Theorem  \ref{powerthm}
is proved in \S\ref{pf13} and
Theorem
\ref{BOext} in \S\ref{sectBOext}.

\medskip

{\it Acknowledgement.} This paper relies substantially on
 ideas in the articles by M. Christ \cite{Ch}, S. Drury \cite{D}
 and by S. Drury and B. Marshall \cite{DM1}, \cite{DM2}. We have added
 an exposition of some of their work hoping that the paper becomes
 more accessible.
We also thank the referees for their comments.

\section{Preliminaries} \label{prelim}

{\bf Lorentz spaces.} We use the standard quasi-norm on the Lorentz space
$L^{p,q}$, namely for $p,q<\infty$
$$\|f\|_{p,q}= \Big(\frac {q}{p}\int_0^\infty \big[ t^{1/p} f^*(t)\big]^q
 \frac{dt}t\Big)^{1/q}, \quad q<\infty
$$
where $f^*$ is the nonincreasing rearrangement of $f$. Moreover
$$\|f\|_{p,\infty}=
\sup_{t>0} t^{1/p} f^*(t)= \sup_{\la>0} \la \big[
\meas\big(\{x: |f(x)|>\la\}\big)\big]^{1/p}.
$$
This does not define  a norm unless
$1\le p=q$; however  $L^{p,q}$ is normable if
$1<p<\infty$ and $1\le q\le \infty$.
For this and many other useful properties on Lorentz spaces we refer
to \cite{hunt} or \cite{st-gw}.

We state some facts on Lorentz spaces needed later.
First, there is the following variant of H\"older's
 inequality for the $L^{r,\infty}$ quasinorms, namely
\begin{equation}\label{hoelderinfty}
\Big\|\prod_{i=1}^n h_i\Big\|_{r,\infty}\le n^{1/r}
\prod_{i=1}^n \big\|h_i\big\|_{rs_i,\infty}, \quad r>0, \,\, \sum_{i=1}^n
\frac{1}{s_i}=1.
\end{equation}
This follows by observing that the set
$\Omega=\{x:\prod_{i=1}^n |h_i|>\la\}$ is contained in the union
of the $n$ sets
$$\Omega_j:=\Big\{x: \frac{|h_j(x)|}{\|h_j\|_{rs_j,\infty}}
>\Big( \frac{\la}{\prod_{i=1}^n \|h_i\|_{r s_i,\infty}}\Big)^{1/s_j}\Big\}.$$

Next, as mentioned in the introduction, we shall use
bounds for  multilinear operators  on $L^{r,\infty}$ for $r<1$.
The advantage of working with the  spaces $L^{r,\infty}$, $r<1$
(as opposed to  $L^{1,\infty}$, say) is that they are
{\it $r$-convex}, that is,
the inequality
\begin{equation}\label{rconvexity}
\Big\|\sum_{l=1}^N h_l\Big\|_{X} \le C_r
\Big(\sum_{l=1}\| h_l\|^r_{X}\Big)^{1/r}
\end{equation}
holds for $X=L^{r,\infty}$ with $C_r$ independent of $N$.
This is a result  of  Stein, Taibleson and G. Weiss \cite{stw} who prove
\eqref{rconvexity} with $C_r=(\tfrac{2-r}{1-r})^{1/r}$; independently
the $r$-convexity of  $L^{r,\infty}$ was shown by
 Kalton \cite{kalton} (who states that  Pisier and Zinn also proved an equivalent result).
Note that in contrast  $L^{1,\infty}$ is not $1$-convex, however a useful
precursor to \eqref{rconvexity} in this case had been found by Stein and N. Weiss in \cite{st-nw}.

We note the following immediate  consequence of  \eqref{rconvexity}.
\begin{lemma} \label{wtr}
 Let $0<r<1$,   $r\le p<\infty$,
and let $Y$ be a complete  $r$-convex space of measurable functions.
Suppose that the  linear operator $S$
maps simple functions
to measurable functions in $Y$ and that for every
for every  measurable set $E$, and for every simple  $f$ with
$|f(x)|\le \chi_E(x)$
 a.e., the inequality $\|Sf\|_{Y}\le C |E|^{1/p}$ holds.
Then $S$ maps $L^{p,r}$ boundedly to $Y$.
\end{lemma}
In particular we may consider  $Y=L^{r,\infty}$, $r<1$; thus
a linear operator $S$ which is  restricted weak type $r$  is also  of
weak type $r$.
To verify Lemma \ref{wtr} let
$f^*$ be the nonincreasing rearrangement of $f$ and let
$E_l$ be the set of all $x$ for which $f^*(2^{l})<|f(x)|\le f^*(2^{l-1})$;
thus the measure of
$E_l$ does not exceed $2^l$.
Define $g_l(x)=  \chi_{E_l}(x) f(x)/|f^*(2^{l-1})|$
if $E_l$ has positive measure (otherwise put $g_l=0$).
By assumption
$\|Sg_l\|_{Y}\le C|E_l|^{1/p}$. Since $f=\sum_l f^*(2^{l-1}) g_l$
we get by  \eqref{rconvexity}
$\|Sf\|_{Y} \le C(\sum_l[f^*(2^{l-1})]^r |E_l|^{r/p})^{1/r}
\le C (\sum_l[2^{l/p} f^*(2^{l-1})]^r)^{1/r}$
which implies the assertion of the Lemma.

We shall use  an analytic interpolation theorem for
Lorentz spaces, for all parameters; this is due to Y. Sagher \cite{sagher}
who extended a version of the Riesz-Thorin theorem for
Lorentz spaces by Hunt \cite{hunt}. These results were proved for all indices
using the harmonic majorization of subharmonic functions.
We state the following  consequence of Sagher's theorem:
\begin{proposition}\label{saghercor}
Let $T$ be a multilinear map defined on $n$ tuples of simple functions,
with values in measurable functions (on some measure space) so
that the inequality
$$
\big\|T(f_1,\dots, f_n)\big\|_{r_i, s_i}
\le A_i  \prod_{j=1}^n \|f_j\|_{p_{j,i}, q_{j,i}}
$$
holds for $i=0,1$.
Then there is  a constant $C$ depending only on the exponents
$r_i, s_i, p_{j,i}, q_{j,i}$ so that
$$
\big\|T(f_1,\dots, f_n)\big\|_{r, s}
 \le C A_0^{1-\vartheta}A_1^\vartheta  \prod_{j=1}^n \|f_j\|_{p_{j}, q_{j}}
$$
holds for
$\big(\tfrac{1}{p_j}, \tfrac{1}{q_j}, \tfrac{1}{r}, \tfrac{1}{s}\big)=
(1-\vth)
\big(\tfrac{1}{p_{j,0}}, \tfrac{1}{q_{j,0}}, \tfrac{1}{r_0}, \tfrac{1}{s_0}\big)+\vth
\big(\tfrac{1}{p_{j,1}}, \tfrac{1}{q_{j,1}}, \tfrac{1}{r_1},
\tfrac{1}{s_1}\big).$
\end{proposition}

This follows from Sagher's theorem by
normalizing each entry $f_j$ in $L^{p,q}$, so that
$\|f_j\|_{p,q}=1$  and then imbedding
each entry
in an analytic family $f_{j,z}$, so that $f_j=f_{j,\vth}$,
 $f_{j,z}= e^{i\arg(f_j)}
G_{j,0}^{1-z} G_{j,1}^z$
with suitable $G_{j,0}$,
$G_{j,1}$, and
$\|G_{j,i}\|_{p_i,q_i}\le C$, $i=0,1$.
See \cite{hunt} and also \cite{sagher}.

\medskip

We shall use a version of a multilinear interpolation argument introduced by
 M. Christ in  \cite{Ch}, often referred to as the multilinear trick.
The result is summarized in

\begin{proposition}\label{trick} Let $0<\beta_j<\infty$, $j=1,\dots, n$, let
$$H=\{x\in \Bbb R^n: \sum_{j=1}^n \beta_j^{-1} x_j =1\},$$  and
let $K$ be a compact subset of $H\cap [0,\infty)^n$. Denote by
$\conv(K)$ the convex hull of $K$ and by
$(conv K)^o$ its interior (with respect to the subspace topology on $H$ induced by $\bbR^n$).

Let $r\le 1$, let $p_j\ge r$, $j=1,\dots, n$, and  let $Y$ be an
$r$-convex Lorentz space (i.e. if $r<1$ then $Y=L^{\rho, q}$ for
$\rho\in [r,\infty)$, $q\in [r,\infty]$, or $Y=L^\infty$). Let  $T$  be  an $n$-linear map
with values in $Y$, defined on $n$-tuples of simple functions,
so that
\begin{equation}\label{pjrassumption}
\|T(f_1,\dots, f_n)\|_Y\le \prod_{j=1}^n \|f_j\|_{L^{p_j,r}},\quad
(p_1^{-1},\dots, p_n^{-1})\in K,
\end{equation}
here, $p_j=\infty$ is allowed but $L^{\infty,r}$ should
be interpreted as $L^\infty$.

Then, for $(p_1^{-1},\dots, p_n^{-1}) \in \big(\conv (K)\big)^o$ and
$ \sum_{j=1}^n q_j^{-1} =r^{-1}$,
\begin{equation}\label{mintconcl}
\|T(f_1,\dots, f_n)\|_Y \le C \prod_{j=1}^n \|f_j\|_{p_j,q_j}.
\end{equation}

\end{proposition}

For Banach spaces $Y$ this is due to Christ \cite{Ch}.
The version for all Lorentz spaces can be
proved using  results on analytic interpolation in the form of
Proposition \ref{saghercor}  in combination with Christ's method and
Lemma \ref{wtr}.
We sketch Christ's argument for the case $r<1$.
First, Proposition \ref{saghercor} yields the inequality
\eqref{pjrassumption} for all $(p_1^{-1},\dots, p_n^{-1})\in \conv(K)$.
Next, the main idea is to assume that $|f_1(x)|\le \chi_E(x)$ a.e.,
and then to prove, for $k=1,\dots, n$, and all
$(p_1^{-1},\dots, p_n^{-1})$ in the interior of $\conv(K)$,
the inequality
\begin{equation}\label{middleint}
\|T(f_1,\dots, f_n)\|_Y \le C(p_1,\dots, p_n) |E|^{1/p_1}
\prod_{j=2}^k \|f_j\|_{p_j,\infty}
\prod_{j=k+1}^n \|f_j\|_{p_j,r};
\end{equation}
with the obvious interpretation that the first product is $1$ if $k=1$ and the second product is $1$ if $k=n$ (we are interested in this last case).
We argue by induction and assume that
\eqref{middleint} is true for some $k= k_o<n$ (the case $k=1$  has been already obtained, in all of $\conv(K)$).
We freeze $p_i$ for $i\notin \{1,k_o+1\}$ and
consider the line segment $\ell$
obtained by intersecting $\conv(K)$ with the two-dimensional plane
$\{x: x_i=p_i^{-1}, i\neq 1, i\neq k_o+1\}$. If $X=(p_1^{-1},\dots,p_d^{-1})$ is in the interior of $\conv(K)$ then it is in the interior of that line segment.
We interpret the inequality
\eqref{middleint} as a linear operator acting on $f_k$ and, by
real interpolation ({\it i.e.} the Marcinkiewicz theorem in its general form)
we get \eqref{middleint}  for $k=k_o+1$ on the open line segment.

Let $\fS^n$ be the group of permutation  on $n$ letters.
Given any  $\varpi\in \fS^n$  we can apply \eqref{middleint} for $k=n$  to the operator $T^\varpi$ defined by
$T^{\varpi}(f_1,\dots, f_n)=
T(f_{\varpi(1)},\dots, f_{\varpi(n)})$,  with $K$  modified appropriately.
By using also  Lemma \ref{wtr} we get,
for $k=1,\dots, n$,
\begin{equation*}
\|T(f_1,\dots, f_n)\|_Y \le C'(p_1,\dots, p_n)
\|f_k\|_{p_k, r}\prod_{j\neq k}  \|f_j\|_{p_j,\infty}.
\end{equation*}
This is already a special case of the assertion and the general case
follows by further  multiple applications of Proposition \ref{saghercor}.

{\it Remark:} Alternatively a more general result can be obtained for
Lions-Peetre interpolation spaces; an elegant  version for  $r$-convex
quasi-normed spaces which in several respects is more general
 is due to  Janson
\cite{janson}, and Proposition \ref{trick} can be seen as  a special case
of his result.

\medskip

{\bf Vandermonde operators.} We now discuss a result
by Drury and Marshall which  concerns certain multilinear
operators involving
the Vandermonde determinants. For a  vector $x\in \Bbb R^d$ let
$V_d(x)$ be the determinant of the
$d\times d$ Vandermonde matrix $(x_i^{j-1})_{i,j}$; {\it i.e.}
\begin{equation} \label{vanderm}
 V_d(x)= \prod_{1\le i<j\le d} (x_j-x_i).
\end{equation}
For $h=(h_1,\dots, h_{d-1})\in (\Bbb R_+)^{d-1}$ define
$\kappa(h)\in [0,\infty)^d$ by
\begin{equation}\label{kappah}
\kappa_1(h)=0, \quad \kappa_{j}(h)= h_1+\dots+ h_{j-1}, \quad 2\le j\le d,
\end{equation}
and
\begin{equation}\label{vh}
v(h)\equiv v_d(h)=V_d(\ka(h)).
\end{equation}

 Define
$$\frak V[f_1,\dots, f_d](t,h)
:=
v(h)^{-1} \prod_{i=1}^d f_i(t+\kappa_i(h)).
$$
Let $L_v^A(L^B)$
denote the weighted mixed norm space
consisting of func\-tions $(t,h)\mapsto G(t,h)$ with
$\|G\|_{L_v^A(L^B)}=(\int \|G(\cdot,h)\|_B^A v(h)dh)^{1/A}<\infty$; then
\begin{multline} \label{mixednormforV}
\|\frak V[f_1,\dots, f_d]\|_{L_v^A(L^B)}=\\
\Big(\int \Big(\int \prod |f_i(t+\kappa_i(h))|^B dt\Big)^{A/B}
v(h)^{1-A} dh\Big)^{1/A}.
\end{multline}

 \begin{proposition}\label{vandthm}  {\it (cf. \cite{DM1}, \cite{DM2})}.

(i) Let, for $\alpha>0$,
$$\Omega_{d}(\alpha)=\{h\in (0,\infty)^{d-1}: v_d(h)\le \alpha \}$$  and assume
$d\ge 2$.
Then  $\Omega_{d}(\alpha)$ has
$(d-1)$-Lebesgue measure $\le C_{d} \alpha^{2/d}$.

(ii)
{\it Suppose that  $1<A<\frac{d+2}d$,  $1<A\le B<\frac{2A}{d+2-dA}$, and
set $\sigma =2/(d+2-dA)$.
For $\nu=1,\dots, d$ let $Q_\nu$ be the point in $\bbR^d$ for which the
$\nu^{\text {th}}$ coordinate is $B^{-1}$ and the other coordinates
are equal to  $(\sigma A)^{-1}$, and let  $\Sigma(A,B)$ be the $d-1$
dimensional closed convex hull of the points $Q_1,\dots, Q_d$. }
{\it Suppose that $(p_1^{-1}, \dots, p_d^{-1})\in \Sigma(A,B)$. Then }
\begin{equation}\label{vandineq}
\big\|\frak V[f_1,\dots, f_d]\big\|_{L_v^A(L^B)}
\le C\prod_{i=1}^d \|f_i\|_{L^{p_i,1}}.
\end{equation}
\end{proposition}

The proof is given in Appendix  \S\ref{app2}.
As has been pointed out in \cite{DM2} the Lorentz spaces
$L^{p_j, 1}$ can be replaced by larger $L^{p_j,q_j}$, provided that
$(p_1^{-1},\dots, p_d^{-1})$ belongs to the interior of $\Sigma(A,B)$
and $\sum_{j=1}^d q_j^{-1}=1$; this follows from Proposition \ref{trick}.
However this improvement of Proposition \ref{vandthm} does not seem to
be relevant
 for the critical
estimates on the extension operator.

\section{Proof of Theorem \ref{osc}: the weak type estimate
}\label{lorentzbd}

Instead of a single oscillatory integral operator it will be convenient
to consider  classes of operators with certain uniform estimates, depending on the derivatives of phase and amplitude.

\medskip

{\bf Definition.}
{\it (i)  Let $N\gg d$ be fixed.
Let $B\ge 3$, and $0<b\le 1/2$.
Denote by $\fA(B)$ the class of functions  $a\in C^{N}(\Bbb R^d\times \bbR)$
which are supported in the cylinder  $\cZ:=\{(x,t): |x|\le 1, |t|\le 1\} $
and
which  satisfy the inequalities
\[
\big|\partial_x^{\alpha} \partial_t^j a(x,t)|\le B, \,\,
|\alpha|\le N, j\le N,
\]
for all $(x,t)\in \cZ$.

(ii)
Let  $\cZ_2:=\{(x,t): |x|\le 2, |t|\le 2\} $
and let $\Phi[B,b]$  be the class of phase functions
$\phi\in C^{N}$
for which the inequalities
$$
|\det\big(  \partial_{t}(\nabla_{x}\phi)  ,\,\partial_{t}^{2}
(\nabla_{x}\phi)  ,\,\ldots,\,\partial_{t}^{n}(\nabla_{x}\phi)  \big)
|\ge b$$
and
\[
\big|\partial_x^{\alpha} \partial_t^j \phi(x,t)|\le B, \quad
1\le j\le N, \quad 1\le |\alpha|\le N,
\]
hold, for all $(x,t)\in \cZ_2$.

(iii)  Let
\begin{equation}
\label{Adef}
\cA_R(B,b):= \sup_{\la\le R} \sup
\big\{ (1+\la)^{d/q_d}
\|T_\la f\|_{q_d,\infty}\big\}
\end{equation}
where the inner supremum is taken over all
$f\in L^{q_d}$ with  $\|f\|_{L^{q_d}}\le 1$  and all oscillatory
integral operators $T_\la$ of the form \eqref{Tla}
for which
the amplitude
$ a$ belongs to $\fA[B]$  and the phase
$\phi$ belongs to  $\Phi[B,b]$.
}

\medskip

Clearly
$\cA_R(B,b)$ is increasing in $R$ and  finite for any choice of $R,B,b$;
an immediate estimate is $\cA_R(B,b)=O(R^{d/q})$ as $R\to\infty$.
However we  need to prove that
$$\cA_R(B,b)=O(1), \quad R\to\infty,$$
with the implicit constant only depending on $B$, $b$.
By suitable rescaling arguments
one can show some stability of the constants $\cA_R(B,b)$, namely:
\begin{lemma}\label{stab} The inequality
$$\cA_R(B,b)\lc C(B,b)  (1+\cA_R(3,2^{-1}))$$
holds for $b\le 1/2$, $B\ge 3$.
\end{lemma}
The notation $\lc$ indicates a constant which  does not depend on
the parameters $B,b, R$ (but may depend on the dimension). It is easy
to see that
Lemma \ref{stab}  is equivalent with
\begin{equation}\label{equivstab}
\cA_R(2B,2b)\lc \widetilde C(B,b)  (1+\cA_R(B,b))
\end{equation}
for $b\le 1/2$, $B\ge 3$.
We shall first
 take Lemma \ref{stab} for granted and give a proof at the end of this section.

\begin{proof}[Proof of Theorem \ref{osc}]
We need to  show that  $\cA_R(B,b)$ in \eqref{Adef} is
bounded, uniformly  in $R$.
We may  assume that the support of $a$ is in a small ball of radius
$\le (b/B^d)^{1000}$. We fix $R\gg1$ and we shall estimate
the $L^{q_d}\to L^{q_d,\infty}$ operator norm
$\|T_\la\|$ for $\la\le R$, provided that $a\in \fA(B)$ and
 $\phi\in \Phi(B,b)$.
We may assume that $\la\ge C_0(B,b)$ (suitably chosen).

We now choose an integer $n>q_d$ and estimate the $n$-linear expression
\begin{equation} \label{fMn}
\fM_\la(g_1,\dots, g_n)= \prod_{j=1}^n T_\la g_j;
\end{equation}
 observe that
$\|T_\la f\|_{q_d,\infty}= \|\fM_\la(f,\dots, f)\|_{q_d/n,\infty}^{1/n}$.
We will take  $n=d q_d$ (although a similar argument works for any $n>q_d$).
Using the symmetry of $\fM_\la$ we may factor it in various ways and
first derive estimates for the
$d$-linear expression
$$\cM_\la(f_1,\dots, f_d)= \prod_{i=1}^d T_\la f_i.$$
Let $\chi_k$ be the indicator function of
\[S_k=\{ t\in \bbR^d :2^{-k-1}\le \prod_{1\le i<j\le d} |t_i-t_j|<2^{-k}
\}.
\]
Following  \cite{baklee} split $\cM_\la=\sum_{k\in \bbZ} \cM_{\la,k}$, where
\[
\cM_{\la,k}(f_1,\dots, f_d)=
\int e^{i\la(\phi(x,t_1)+\dots\phi(x,t_d))}\prod_{i=1}^d \big[
a(x,t_i) f_i(t_i)\chi_k(t)\big] dt_1\dots dt_d.
\]
We first use a by now standard $L^2$ estimate
(for the complete details see \cite{baklee}).
One may introduce the singular change of variables
$y=\sum_{i=1}^d\phi(x,t_i)$, apply a standard $L^2$ estimate for singular integrals and change variables back  and it
follows that
\begin{multline*}
\|\cM_{\la,k}(f_1,\dots, f_d)\|_2\le
\\ C_1(B,b)
 (1+\la)^{-d/2}
\Big(\int_{S_k} 2^k \big|f_1(t_1)\cdots f_d(t_d)|^2 dt_1\cdots dt_d\Big)^{1/2}.
\end{multline*}
The factor $2^k$ is comparable to the reciprocal of
the Jacobian, which is
$(\prod_{1\le i<j\le d} |t_i-t_j|)^{-1}$.
The measure of the $t_1$-section of $S_k$
({\it  i.e.} the set of
 all $(t_2,\dots,t_d)$ for which $(t_1,\dots, t_d)\in S_k$)  is
$O(2^{-2k/d})$ (\cf. Proposition \ref{vandthm}, (i)). Thus also
\begin{equation}\label{L2}
\|\cM_{\la,k}(f_1,\dots, f_d)\|_2\lc C_1(B,b) 2^{(d-2)k/(2d)}
 (1+\la)^{-d/2}
\|f_1\|_2 \prod_{i=2}^d\|f_i\|_\infty.
\end{equation}

Now let $u(h_1,\dots, h_{d-1})= h_1\cdots h_{d-1}
\prod_{1\le i<j\le d-1} |h_i-h_j|$.
Then
$$
\cM_{\la,k}(f_1,\dots, f_d)=
\int_{h: 2^{-k-1}\le u(h)\le 2^{-k}}
\int e^{i\la \psi(x,s, h)} \fa(x,s,h)
F(s,h)
ds dh
$$
where
$$F(s,h)=f_1(s)\prod_{i=1}^{d-1} f_{i+1}(s+h_i),$$
$$\psi(x,s,h) = \phi(x,s)+ \sum_{i=1}^{d-1}\phi(x,s+h_i),$$ and
$$\fa(x,s,h)=a(x,s)\prod_{i=1}^{d-1} a(x,s+h_i).$$
Then it is easy to check that
$$
\fa (\cdot,h) \in \fA(  CB),\quad
\psi(\cdot,h)\in \Phi[CB,C^{-1}b],
$$
for some absolute constant $C$, uniformly for the
$h$ under consideration.

Define
$$S_{h,\la} F(x)=\int e^{i\la \psi(x,t,h)}\fa(x,t,h) F(t) dt$$
then by the definition of $\cA_R(CB, b/C)$
$$(1+\la)^{d/q_d}\|S_{h,\la} F\|_{q_d,\infty}\le
\cA_R(CB, b/C)
\|F\|_{L^{q_d}}.
$$
Using the (integral form of the) triangle inequality
\begin{align*}
(1+\la)^{d/q_d}&\big\|\cM_{\la,k}(f_1,\dots, f_d)\big\|_{q_d,\infty}
\\&\le
\int_{h: 2^{-k-1}\le u(h)\le 2^{-k}}
(1+\la)^{d/q_d}\|S_{h,\la} F(\cdot,h)\|_{q_d,\infty} dh
\notag
\\ &\le
\cA_R(CB, b/C)
\int_{h: 2^{-k-1}\le u(h)\le 2^{-k}}
\|F(\cdot,h)\|_{q_d}dh
\notag
\\ &\le \cA_R(CB, b/C)
\int_{h: 2^{-k-1}\le u(h)\le 2^{-k}}
\|f_1\|_{q_d} \prod_{i=2}^d\|f_i\|_\infty \,dh
\notag
\\ &\le
C' \cA_R(CB, b/C) 2^{-2k/d}
\|f_1\|_{q_d} \prod_{i=2}^d\|f_i\|_\infty;
\end{align*}
here we used again, that
$\meas(\{h: |u(h)|\le  2^{-k}\})\approx
2^{-2k/d}$.

By Lemma
\ref{stab} (\cf. \eqref{equivstab}) we also get
\begin{multline}\label{qdendpoint}
(1+\la)^{d/q_d}\big\|\cM_{\la,k}(f_1,\dots, f_d)\big\|_{q_d,\infty}
\\
\le
C_2(B,b)  (1+\cA_R(B, b)) 2^{-2k/d}
\|f_1\|_{q_d} \prod_{i=2}^d\|f_i\|_\infty.
\end{multline}

Now we interpolate the $L^2$ and $L^{q_d,\infty}$ bounds \eqref{L2} and
\eqref{qdendpoint}
by the
real method (based on the familiar  argument by Bourgain
\cite{Bo1}  for the
spherical maximal function, see also \cite{CSWW}).
Let $\vth\in (0,1)$ be defined by
$$(1-\vth)
\frac{d-2}{2d} +\vth (\frac{-2}d)=0,$$ then
notice that
\begin{equation}\label{thetarelation}
\vth=\frac{d-2}{d+2}
\text{ and }
 \frac{1-\vth}{2}+\frac{\vth}{q_d}= \frac{d}{q_d},
\end{equation}
and $0<\vth<1$ if $d\ge 3$.

Thus, for fixed $f_2, \dots, f_d$, the linear operator
 $f_1\mapsto
\sum_k\cM_{\la,k}(f_1,\dots, f_d)$ is of restricted weak type
$(q_d/d, q_d/d)$; in fact
\begin{align*}
&\Big\|
\sum_k \cM_{\la,k}(f_1,\dots, f_d)\Big\|_{q_d/d,\infty}
\le C_3 [C_1(B,b) (1+\la)^{-d/2}]^{1-\vth}
\\
&\quad\quad\quad\times
[(1+\la)^{-d/q_d}(1+\cA_R( B, b))]^\vth
\|f_1\|_{q_d/d, 1} \prod_{i=2}^d\|f_i\|_\infty
\\&\le
C_4(B,b) (1+\la)^{-d^2/q_d} (1+\cA_R(B, b))^\vth
\|f_1\|_{q_d/d, 1} \prod_{i=2}^d\|f_i\|_\infty.
\end{align*}
By the symmetry and various interpolations this estimate
leads to a restricted weak type
(or even improved Lorentz type) estimate for $T_\la$; however to  prove the
stronger weak type
estimate we now set $n=d q_d$ and  consider the $n$-linear  operator
\eqref{fMn}.



We use \eqref{hoelderinfty}
to bound
\begin{multline*}
\Big\|\prod_{i=1}^{dq_d} T_\la g_i\Big\|_{1/d,\infty}
\le C \prod_{k=0}^{q_d-1}
\Big\|
 \cM_{\la}(g_{kd+1},\dots, g_{(k+1)d})
\Big\|_{q_d/d,\infty}
\\
\le C_4(B,b)^{q_d} (1+\la)^{-d^2} (1+\cA_R(B, b))^{\vth q_d}
\prod_{k=0}^{q_d-1} \Big[\|g_{kd+1}\|_{q_d/d, 1}
\prod_{i=2}^d\|g_{kd+i}\|_\infty\Big].
\end{multline*}

Using the symmetry of the operator
we get for any permutation $\pi$ on $n=d q_d$ letters
\begin{align}\label{symmetrization}
\Big\|\prod_{i=1}^{dq_d} T_\la g_i\Big\|_{1/d,\infty}
\le C_5(B,b) (1+\cA_R(B,b))^{\vth q_d}(1+\la)^{-d^2}
\prod_{j=1}^{dq_d} \|g_{\pi(j)}\|_{p_j,r_j}
\end{align}
where $(p_1^{-1}, \dots, p_{dq_d}^{-1})$ is in the set
$$K=\{
P^\pi=(P^\pi_1,\dots, P^\pi_{q_d}): \pi \in \fS_{dq_d}\},$$ with $P^\pi$  defined by
\[
P^\pi_{\pi(i)}=\begin{cases}
d/q_d, &i=1,\dots, q_d,
\\
0,\quad &q_d+1\le i\le dq_d,
\end{cases}
\]
and with $r_j=1$ if $1\le j\le q_d$ and $r_j=\infty$, $j>q_d$
 (so that $L^{\infty,\infty}=L^\infty$ in \eqref{symmetrization}).
The convex hull of the  set $K$ is a simplex
on the hyperplane
$\{X\in \bbR^{dq_d}:\sum_{i=1}^{dq_d} X_i=d\}$ with vertices $P^{\pi}$.
Using the multilinear interpolation result of Proposition \ref{trick},
part (ii), we obtain that
\begin{align}\label{fullafterinterpol}
\Big\|\prod_{i=1}^{dq_d} T_\la g_i\Big\|_{1/d,\infty}
\le C_5(B,b) (1+\cA_R(B,b))^{\vth q_d}(1+\la)^{-d^2}
\prod_{j=1}^{dq_d} \|g_{\pi(j)}\|_{p_j,s_j}
\end{align}
for all
$(p_1^{-1}, \dots, p_{dq_d}^{-1}) \in (\conv(K))^o$ and
$\sum_{j=1}^{dq_d}s_j^{-1} = r^{-1}\equiv d$.
The center $(q_d^{-1},\dots, q_d^{-1})$ lies  in
$(\conv(K))^o$ and we can choose $s_j=q_d$, for $j=1,\dots, dq_d$.

Setting $g_j=f$ we get
\begin{align*}
\big\|T_\la f\big\|_{q_d,\infty}^{dq_d}&=
\big\|[T_\la f]^{d q_d}\big\|_{1/d,\infty}
\\&\le C_5(B,b) (1+\cA_R(B,b))^{\vth q_d}(1+\la)^{-d^2}
\|f\|_{q_d}^{dq_d}.
\end{align*}
Thus, by definition of $\cA_R$,
\[ \cA_R(B,b)\le
C_6(B,b)
(1+\cA_R(B,b))^{\vth/d}
\]
which gives
$\cA_R(B,b)=O(1)$ as $R\to \infty$.
\end{proof}

\medskip

\noi{\it Remark:} It is perhaps instructive to  compare this
argument  with one in differentiation theory, namely
Christ's simplification of the $L^p$ boundedness result ($p>1$) by
Nagel, Stein and Wainger \cite{na-s-w} on differentiation in
lacunary directions. In \cite{na-s-w} a bootstrap argument
is used which is formally similar to Drury's argument.
Our argument resembles the simplification   which avoids this
iteration, see {\it e.g.} Theorem B in  \cite{Ca} for an exposition.
\medskip

\begin{proof}[Proof of Lemma \ref{stab}]
Let $\phi$ in $\Phi(B,b)$ and  let  $a\in \fA(B)$. We wish to bound
the $L^{q}\to L^{q,\infty}$ operator norm of $T_\la$ where
$$T_\la f(w)= \int e^{i\la \phi(w,s)} a(w,s) f(s) ds.$$

Let $\chi\in C^\infty(\bbR)$ be a function supported in $(-1,1)$
satisfying $$\sum_{n\in \bbZ}\chi(s-n)=1$$ for all $s\in \bbR$.
The argument will involve rescalings depending on two
small numbers $\eps \ll 1$ and $ \delta \ll\eps$, in fact
we shall see that
\begin{align}\label{defeps}
\eps&=  \big( 10 (d!)^2 B^d b^{-1}\big)^{-d},
\\
\label{defdel}
\delta&= 10^{-2Nd^2} (d!)^{-1}
(d!B^d b^{-1})^{-Nd^2}
\end{align}
is an admissible choice.

We begin by observing the trivial estimate
\begin{multline} \label{trivial}
\|T_\la\|_{L^{q}\to L^{q,\infty}}
\le C_d \le C_d'(1+\delta^{-1})^{d/q} (1+\la)^{-d/q}, \\
\quad \text{ if } \la \le \min\{\delta^{-1}, R\}.
\end{multline}
This takes care of the case $\la\le \delta^{-1}$,
and in what follows we shall assume
 that $\la\ge \delta^{-1}$. We decompose the symbol into pieces supported
 in $(\tfrac{\delta} 2,\dots, \tfrac{\delta} 2, \eps)$ boxes. To this end we set, for $\mu\in \bbZ$ and $\nu\in \bbZ^d$,
$$a_{\mu\nu}(w,s)= a(w,s)
\chi_d (4\delta^{-1} w-\nu)\chi(2\eps^{-1}s-\mu)
$$ where
$\chi_d(w):=
\prod_{i=1}^d
\chi( w_i)$. We also set  $s(\mu)= \eps\mu/2$, $w(\nu)=\delta \nu/4$, $P_{\mu\nu}= (w(\nu),s(\mu)).$
Let $J^{\mu\nu}_\phi$ be the $d\times d$ matrix with
$$\big(J^{\mu\nu}_\phi\big)_{jk}= \partial_s^j \partial_{w_k}\phi(P_{\mu,\nu}).$$
By assumption $|\det J^{\mu\nu}_\phi|\ge b$.
Let $L^{\mu\nu}$ be the inverse matrix of
$J^{\mu\nu}_\phi$. By the cofactor formula we see that its entries
have the bound
\begin{equation}\label{Lbd}
\big|[L^{\mu\nu}]_{jk}\big|\le (d-1)! B^{d-1} b^{-1}.
\end{equation}
We also denote by
$J^{\mu\nu}_\phi[\eps]$
the matrix whose $i^{\text{th}}$ row  is
$\eps^i \partial_s^i \nabla_w \phi(P_{\mu\nu})$. Let
$L^{\mu\nu}[\eps]$
be its inverse so that  the $k^{\text{th}}$ column of
$L^{\mu\nu}[\eps]$  is equal to $\eps^{-k}$ times the
$k^{\text{th}}$ column of $L^{\mu\nu}$.
Let
$ \cT^{\mu\nu}_\la$ be the oscillatory integral operator with
phase $\phi$ and amplitude  $a_{\mu\nu}$.
Then, by the support properties of $a$ and $a_{\mu\nu}$,
\begin{equation} \label{summing}\|T_\la f\|_{q,\infty} \le 10^{d+1}
\delta^{-d} \eps^{-1} \max_{\mu,\nu} \|\cT_\la^{\mu\nu} f\|_{L^{q,\infty}}
\end{equation}
 and it suffices to estimate the individual operators $\cT_\la^{\mu\nu}$.

In what follows we fix $\mu,\nu$ and usually
drop the superscripts $\mu,\nu$ in $L^{\mu\nu}$ and $L^{\mu\nu}[\eps]$.

Define a rescaled operator $S_\Lambda\equiv S^{\mu\nu}_\La$ by
$$S_\La g(x) =
\int e^{i\La \psi(x,t)} u(x,t) g(t) dt
$$
where
\begin{equation}
\psi(x,t) = \delta^{-1} \phi( w(\nu)+ L[\eps] \delta x,
s(\mu)+ \eps t)
\label{psidef}
\end{equation}
and
\begin{align} \label{udef}
u(x,t)&\equiv a_{\mu\nu} ( w(\nu)+\delta L[\eps]  x, s(\mu)+ \eps t)
\\&=
a( w(\nu)+ \delta L[\eps] x, s(\mu)+ \eps t)
\chi(2t)
\chi_d(4L[\eps] x) .
\notag
\end{align}

Then
$$\cT_\la^{\mu\nu}  f(w(\nu)+  \delta L[\eps]  x)= \eps S_{\la \delta} [f(s(\mu)+
 \eps \cdot)] (x)
$$
and it follows
\be{Tmunula}\|\cT_\la^{\mu\nu}\|_{L^{q}\to L^{q,\infty}}
 \le \eps^{1-1/q}  |\det L[\eps]|^{1/q} \delta^{d/q}
\|S_{\la\delta} \|_{L^{q}\to L^{q,\infty}}.
\ee

We verify that
$|\det L[\eps]| \le b^{-1} \eps^{-d(d+1)/2}$, by the lower bound for
$\det J^{\mu\nu}_\phi$ and that
$\|L[\eps]\|_\infty  \le \eps^{-d}(d-1)! B^{d-1} b^{-1}$; here
 $\|L[\eps]\|_\infty:=\max_{i,j}|L_{ij}[\eps]|$.
 We shall then check that
\begin{equation} \label{normalized}
\psi \in \Phi(3,1/2)
\end{equation}
and
\begin{equation} \label{normalized2}
\frac{u}{
 \|\chi\|_{C^N}^{d+1} (8d! d B^d b^{-1}\eps^{-d})^N }
\in \fA(1).
\end{equation}
Given \eqref{normalized} and \eqref{normalized2} it follows that,
 for $\delta^{-1} \le \la\le R$,
$$\|S_{\la\delta} \|_{L^{q}\to L^{q,\infty}}\le
(1+\la\delta)^{-d/q}
 \|\chi\|_{C_{N}}^{d+1} (8d! d B^d b^{-1}\eps^{-d})^N
\cA_R(3,1/2)
$$ and thus, combining this estimate with \eqref{summing} and
\eqref{Tmunula} we obtain
\begin{multline*}
(1+\la)^{d/q}\|T_\la\|_{L^{q}\to L^{q,\infty}}
\le \delta^{-d/q} (1+\la\delta)^{d/q}\|T_\la\|_{L^{q}\to L^{q,\infty}}
\\\le  10^{d+1}\delta^{-d} \eps^{-(d^2+d+2)/2q} b^{-1/q}
 \|\chi\|_{C^N}^{d+1} (8d! d B^d b^{-1}\eps^{-d})^N
\cA_R(3,1/2).
\end{multline*}
for $\la\ge \delta^{-1}$.
If we also take into account the trivial bound \eqref{trivial}, and the dependence of $\delta$ and $\eps$ on $B$ and $b$  then
we get
$$\cA_R(B,b)\le C(B,b,N,d) \cA_R(3,1/2).$$

It remains to check \eqref{normalized} and
\eqref{normalized2}. The latter follows by straightforward applications of the
Leibniz rule. Concerning \eqref{normalized}
we consider the matrix
$J_\psi$ with entries $[J_\psi(x,t)]_{jk} =\partial_t^{j} \psi_{x_k}(x,t) $. By definition
$J_\psi(0,0)$ is  the identity matrix.
We expand using Taylor's formula
\begin{equation}\label{taylexp}
\partial_t^j \psi_{x_k}(0,t)=\sum_{l=0}^{d-j} \frac{t^l}{l!}
\partial_t^{j+l}\psi_{x_k}(0,0)
+ \frac{t^{d-j+1}}{(d-j+1)!}
\partial_t^{d+1}\psi_{x_k}(0,\tilde t)
\end{equation}
and notice that
the first sum equals (with $P\equiv P_{\mu\nu}$)
\[
\sum_{l=0}^{d-j} \frac{t^l}{l!} \sum_{m=1}^d \partial_s^{j+l}\phi_{w_m}(P)\eps^{j+l}
L_{mk}\eps^{-k}=
\begin{cases} 0 & \text{if } k<j
\\
\tfrac{t^{k-j}}{(k-j)!}\,
&\text{if } k\ge j.
\end{cases}
\]
The absolute value of the remainder term in \eqref{taylexp}  is
\[
\Big|
\frac{t^{d-j+1}}{(d-j+1)!}\sum_{m=1}^d \partial_s^{d+1}\phi( \tilde w,\tilde s) \eps^{d+1} L_{mk} \eps^{-k}\Big|
\le 2dB\eps \|L\|_\infty \le 2d! B^d b^{-1}\eps.
\]
There is also another error term for the expansion  about $x=0$, and
we have
\[
\partial_t^j\psi_{x_k}(x,t)- \partial_t^j\psi_{x_k}(0,t)
=\sum_{l=1}^d \partial_t^j \psi_{x_kx_l} (\tx, t ) x_l
\]
with
\begin{align*}
\big|\partial_t^j \psi_{x_kx_l}
(\tx, t) \big|&= \big|\delta^{-1} \eps^{j-k-l}\delta^2  \sum_{m,n=1}^d
L_{ml}
L_{nk}
\partial_s^j\phi_{w_nw_m}(\tw, \ts) \big|
\\
&\le \delta^{-1+2} \eps^{j-l-k} Bd^2  (b^{-1}(d-1)! B^{d-1})^2 \le
\delta \eps^{j-l-k} (d!B^d b^{-1})^2.
\end{align*}
Thus, for all $(x,t)\in \cZ_2$,
\begin{equation}\label{Jpsi}
\big[J_\psi(x,t)\big]_{jk}=
\begin{cases}
\frac{t^{k-j}}{(k-j)!} +E_{jk}(x,t), \ &\text{if } k\ge j
\\
E_{jk}(x,t), \ &\text{if } k< j,
\end{cases}
\end{equation}
with
\begin{equation}\label{Ejk}
|E_{jk}(x,t)|\le 2d! B^d b^{-1}\eps  +
2\delta\eps^{1-2d} B (b^{-1}d! B^{d-1})^2.
\end{equation}
By straightforward considerations using cofactor expansions we see that
$$|\det J_\psi(x,t)-1|\le (d!-1)\max_{1\le \ka\le d}\max_{jk}
|E_{jk}(x,t)|^\ka,
$$
and thus, by our choice of $\eps$,
we have
\be{detbd}\det J_\psi(x,t) \ge 1/2, \quad (x,t)\in \cZ_2,
\ee
moreover, using also our choice of $\delta$ and the assumption $|t|\le 2$
\be{upperbdderfirst}
\big|\partial_t^j\psi_{x_k}(x,t)\big|\le 3,\quad  1\le j\le N.
\ee

The above estimates for the second derivatives can be extended
in a straightforward manner to higher derivatives and we obtain for
$j\le N$, and multiindices $\alpha=(\alpha_1,\dots, \alpha_d)$
with $|\alpha|:=\sum_{i=1}^d \alpha_i$  that
\begin{equation*}
\big|\partial_t^j\partial_x^\alpha\psi(x,t)\big|\le
\frac {1}{\delta} \eps^{j-\alpha_1-2\alpha_2-\dots-d\alpha_d}
\delta^{|\alpha|} ( d! B^{d} b^{-1})^{|\alpha|}.
\end{equation*}
Observe that when we have at least two $x$-differentiations then
the smallness of $\delta$ can be used. By our choice \eqref{defdel}
\be{upperbdder}
\big|\partial_t^j\partial_x^\alpha\psi(x,t)\big|\le 2,\quad  1\le j\le N,
\,\, 2\le |\alpha|\le N,
\ee
and it follows from \eqref{detbd}, \eqref{upperbdderfirst}
and \eqref{upperbdder} that $\psi\in \Phi[3,2^{-1}]$.
\end{proof}

\noi{\bf The extension operators for nondegenerate curves.}
The model case for our class of phase functions  is
$\phi(x,t)=-\inn{x}{\gamma(t)}$
where
$\gamma:I\to \bbR^d$ is defined on the compact interval $I$ and
has the property that for each $t$ the derivatives
$\gamma'(t)$, ..., $\gamma^{(d)}(t)$ are linearly independent.
Define the Fourier extension operator
$$\cE f(\xi)= \int_I f(t) e^{-i\inn{\xi}{\gamma(t)}} dt$$
and let $B(\la)$ be a ball in $\bbR^d$ of radius $\la$.
Then by a change of variable
Theorem \ref{osc} implies that there is $C>0$ so that for all $\alpha>0$
$$\meas\big(\{\xi\in B(\la):|\cE f(\xi)|>\alpha\}\big) \le
\big[ C\alpha^{-1} \|f\|_{q_d}\big]^{q_d}.
$$
By letting $\la\to \infty$ and using the monotone convergence
theorem we see that $\cE: L^{q_d}(I)\to L^{q_d,\infty}(\bbR^d)$.
A duality argument shows the local version of Theorem \ref{restr},
namely
$$\big\|\widehat f\circ\gamma\big\|_{L^{p_d}(I)}\le c_I
\|f\|_{ L^{p_d,1}(\bbR^d)}.$$
A nonisotropic scaling
using the dilations $x\mapsto (ux_1,u^2 x_2,\dots, u^d x_d)$
can be used to  deduce the  global version of Theorem \ref{restr}.

\section{Proof of the $L^{q_d}$ bound}\label{strongtypebound}

We now  show  \eqref{Lqbd}.
Recall the bounds for $\cM_{\la, k}$.
As $T_\la$ has bounded operator norms the estimate \eqref{L2}
is wasteful for  large $k$ and the term
$2^{k(1/2-1/d)}(1+\la)^{-d/2}$ can be replaced by a constant.

Note that  for all $k\ge 0$
\begin{equation*}
\big\|\cM_{\la,k}(f_1,\dots, f_d)\big\|_{q_d,\infty}\lc (1+\la)^{-d/q_d}
 2^{-2k/d}
\|f_1\|_{q_d,1} \prod_{i=2}^d\|f_i\|_\infty,
\end{equation*}
which follows from
\eqref{qdendpoint} since we have already  established the restricted weak
type bound for $q_d$.
By real interpolation,
\begin{equation}\label{singlek}
\big\|\cM_{\la,k}(f_1,\dots, f_d)\big\|_{q_d/d,1}\lc
(1+\la)^{-d^2/q_d}
\|f_1\|_{q_d,1} \prod_{i=2}^d\|f_i\|_\infty,
\end{equation}
but there is  also the trivial
bound
\begin{equation}\label{trivint}
\big\|\cM_{\la,k}(f_1,\dots, f_d)\big\|_{q_d/d,1}\lc
(1+\la)^{-d\vth/q_d} 2^{-2k\vth/d}
\|f_1\|_{q_d/d,1} \prod_{i=2}^d\|f_i\|_\infty,
\end{equation}
with $\vth=(d-2)/(d+2)$.
Let $$N_\la=10 d^3 \log \la$$ then certainly by \eqref{trivint}
\begin{equation}\label{error}
\Big\|\sum_{k>N_\la} \cM_{\la,k}(f_1,\dots, f_d)\Big\|_{q_d/d,1}\lc
(1+\la)^{-d^2/q_d}
\|f_1\|_{q_d/d,1} \prod_{i=2}^d\|f_i\|_\infty.
\end{equation}

Furthermore one can show, for $1\le \tau\le \infty$,
\begin{equation}\label{kleN}
\Big\|\sum_{0\le k\le N_\la}\cM_{\la,k}(f_1,\dots, f_d)
\Big\|_{q_d/d,\tau}\lc  N_\la^{1/\tau}
(1+\la)^{-d^2/q_d}
\|f_1\|_{q_d/d,1} \prod_{i=2}^d\|f_i\|_\infty.
\end{equation}
This  follows from the case $\tau=1$ which holds by
\eqref{singlek} and the case $\tau=\infty$
which is the restricted weak type estimate that follows from Bourgain's interpolation argument.
All together
\begin{equation} \label{allk}
\Big\| \prod_{j=1}^d T_\la f_j\Big\|_{q_d/d,\tau}\lc  N_\la^{1/\tau}
(1+\la)^{-d^2/q_d}
\|f_1\|_{q_d/d,1} \prod_{i=2}^d\|f_i\|_\infty,
\end{equation}
and similar bounds with the $f_i$ permuted.
We apply this with
$\tau=q_d/d$ and  use the multilinear trick  for the
$q_d$-linear expression
$\prod_{j=1}^{q_d} T_\la f_j$ on $L^1$. We have
 for all permutations $\pi$ on $q_d$ letters
\begin{align*}
\Big\|\prod_{j=1}^{q_d} T_\la f_j\Big\|_{1}
&\le \Big\|\prod_{i=1}^d T_\la f_{\pi(i)}\Big\|_{\frac{q_d}d,1}
\Big\|\prod_{j=d+1}^{q_d} T_\la f_{\pi(j)}\Big\|_{\frac{q_d}{q_d-d},\infty}
\\
&\le \Big\|\prod_{i=1}^d T_\la f_{\pi(i)}\Big\|_{\frac{q_d}d,1}
\prod_{j=d+1}^{q_d} \big\|T_\la f_{\pi(j)}\big\|_{q_d,\infty}
\\&\lc N_\la (1+\la)^{-d}
\|f_{\pi(1)}\|_{q_d/d,1} \prod_{i=2}^d \|f_{\pi(i)}\|_\infty
\prod_{j=d+1}^{q_d} \|f_{\pi(j)}\|_{q_d,1}.
\end{align*}
The multilinear interpolation result of Proposition \ref{trick},
for $Y=L^1$, yields
\begin{equation*}\Big\|\prod_{j=1}^{q_d} T_\la f_j\Big\|_{1} \lc
N_\la (1+\la)^{-d}
 \prod_{i=1}^{q_d} \|f_{\pi(i)}\|_{p_i,r_i}
\end{equation*}
for $(p_1^{-1},\dots, p_{q_d}^{-1})$ in a neighborhood of
$(q_d^{-1},...,q_d^{-1})$,  satisfying $\sum_{i=1}^{q_d}p_i^{-1}=1$,
and for $\sum_{i=1}^{q_d}r_i^{-1}=1$.
Now $N_\la\approx \log \la$ for $\la \ge 2$ and the asserted result follows if we set
$p_i=r_i=q_d$, $f_i=f$.\qed


\section{A lower bound} \label{sharpness}
We show that the extension operator for the nondegenerate case
 does not map  $L^{q_d,r}(I)$ to $L^{q_d,\infty}(\bbR^d)
$ for $r>q_d$.  By the uniform boundedness principle it suffices
to consider smooth and compactly supported functions  and show that the operator norm is not  finite. We may assume that $I=(-1,1)$. By a linear change of variable we may also assume that $\gamma^{(j)}(0)=e_j$, for
$j=1,\dots, d$.

Let $\chi$ be a nonnegative $C^\infty_0$ function supported in
$(-1/8, 1/8)$ with $\chi(t)=1$ for $t\in (-1/10, 1/10)$.
For $n\in \bbN$ define
\begin{equation}\label{fN}
\begin{aligned}
u_n(t)&= 2^{n/q_d} \chi (2^n (t-2^{-n})),
\\
f_N(t)&=\sum_{n={N+1}}^{2N} u_n(t).
\end{aligned}
\end{equation}
It is easy to see that for $N\ge 2$
\begin{equation}\label{fNup}
\big\|f_N\|_{L^{q_d,r}(I)}\le C N^{1/r}\, ,
\end{equation}
and
thus it suffices to show that for large $N$
\begin{equation}\label{fNdown}
\big\|\cE f_N\|_{L^{q_d,\infty}(\bbR^d)}\ge C N^{1/q_d}.
\end{equation}

In order to achieve this we need the following
van der Corput type asymptotics which is taken from Lemma 5.1 in
 \cite{bggist}.

\medskip
{\bf  Asymptotics.}
\textit{
Let $0<h\le 1$, $I=[-h,h]$, $I^*=[-2h,2h]$
 and let  $g\in C^{2}(I^*)$.
Suppose that
 $h\le 10^{-1}(1+\|g\|_{C^2(I^*)})^{-1}$
and
let $\eta\in C^{1}$ be supported in $I$ and satisfy the bounds}
\begin{equation}\label{etaassumption}
\|\eta\|_{\infty}+\|\eta^{\prime}\|_{1} \le A_{0},
\text{ and } \|\eta^{\prime}\|_{\infty}\le A_{1}.
\end{equation}
\textit{Let $k\ge 2$ and define}
\begin{equation}
I_{\lambda}(\eta,x)= \int\eta(s) \exp\big(i\lambda(\sum_{j=1}^{k-2} x_{j}
s^{j}+ s^{k} + g(s) s^{k+1})\big) ds.
\end{equation}
\textit{Let $\alpha_{k}=\tfrac2k\Gamma(\tfrac1k)\sin(\tfrac{(k-1)\pi}{2k}),
$ if $k$ is odd and $\alpha_k=
\tfrac 2k\Gamma(\tfrac1k) \exp(i\tfrac{\pi}{2k}),$ if $k$ is even.
Suppose that
$|x_{j}|\le\eps\lambda^{(j-k)/k}$, $j=1,\dots, k-2$.
Then there is an absolute constant $C$ so that, for $\lambda>2$,}
\[
|I_{\lambda}(\eta,x)-\eta(0) \alpha_{k} \lambda^{-1/k}| \le C[A_{0}
\eps\lambda^{-1/k}+A_{1}\lambda^{-2/k}(1+\delta_{2,k} \log\lambda)];
\]
\textit{here $\delta_{2,2}=1$, and $\delta_{2,k}=0$ for $k>2$.}

\begin{proof}[Proof of \eqref{fNdown}]
We shall  get good  lower bounds for the set
where $|\cE f_N|\ge\beta $ provided that
$\beta  \ll 2^{-2N}$.
Consider large $\xi$ with $\xi_d \approx |\xi|$. By the implicit
function theorem the equation $\inn{\gamma^{(d-1)}(t)}{\xi}=0$
has a unique solution $\tcr(\xi)$ which is homogeneous of degree zero.

For each $n\in [N,2N]$ we let $\la_n= 2^{nd/q_d}\beta^{-d}$ and
set
\begin{multline*}\cV_n=\{\xi: |\xi'|\le c|\xi_d|, \la_n\le |\xi_d|\le 2\la_n, \, \tcr(\xi)\in
(\tfrac 9{10}2^{-n}, \tfrac {11}{10}2^{-n}),
\\
|\inn{\gamma^{(j)}(\tcr(\xi))}{\xi}| \le \eps \la_n^{j/d}, j=1,\dots, d-2
\}.
\end{multline*}

Note that if $t\in \supp (u_k)$ and $k\neq n$ then
$|t-\tcr(\xi)|\ge c2^{-n}$ and therefore
$|\inn{\gamma^{(d-1)}(t)}{\xi}|\ge c' 2^{-n}|\xi|$.
By van der Corput's lemma with $(d-1)$ derivatives
we get the bound
\begin{equation*}
|\cE u_k(\xi)|\le C 2^{n/q_d} 2^{n/(d-1)} \la_n^{-1/(d-1)},
\qquad \xi\in \cV_n, \quad k\neq n.
\end{equation*}

By the asymptotics above, if  $\eps>0$ is sufficiently small,
then
\begin{equation*}
|\cE u_n(\xi)|\ge  c 2^{n/q_d} \la_n^{-1/d}
- 2^{2n /q_d} \la_n^{-2/d}
\ge c' 2^{n/q_d} \la_n^{-1/d}=c'\beta,
\quad \xi\in \cV_n.
\end{equation*}

Combining the last two inequalities
we obtain
\begin{equation}
  |\cE f_N(\xi)|\ge c'' \beta, \qquad\xi \in \cup_{n=N}^{2N}\cV_n.
\end{equation}
The measure of $\cV_n$ is
$\ge c_\eps 2^{-n}\la_n^{\frac{(d-2)(d-1)}{2d} +2}
=c_\eps 2^{-n} \la_n^{q_d/d}=
c_\eps \beta^{-q_d}$
and the sets $\cV_n$, $N\le n\le 2N$ are disjoint if $N$ is large.
Thus for $N$ large
\eqref{fNdown} follows.
\end{proof}

\section{Proof of Theorem \ref{powerthm}} \label{pf13}
We first note that it suffices to assume that the powers $b_i$ are
mutually distinct  and also $b_i\neq 0$; in the other cases the
weight vanishes identically.

 We only need to prove the result for
$I=(0, 1]$ here, by a scaling argument we can easily extend
the result to $I=(0,\infty)$, using the linear isomorphisms
$x\mapsto (s^{b_1}x_1, ...., s^{b_d}x_d)$. Following \cite{DM2} we will use the
exponential parametrization,  replacing $t$ by $e^{-t}$.
Setting $a_i=-b_i$  we
may assume, after a further linear change of variables,  that
\begin{equation}\label{expcurv}
\gamma(t) = (a_1^{-1}e^{a_1 t}, \cdots, a_d^{-1}e^{a_d t}), \quad 0<t< \infty,
\end{equation}
%
where the $a_j$ are real numbers so that either
(i) $a_1<\dots<a_d<0$,  or
(ii) $0<a_1<\dots<a_d$, or
(iii)  $a_1< \cdots< a_m <0 <
a_{m+1} <\cdots <a_d$, for some $m \in \{ 1, \cdots, d\}$.
We shall give the argument for case (iii), and the proofs for the
other cases require
only notational changes.

Fix any point $(1/p, 1/q)$ on the critical line segment
$1/p+(d^2+d)/(2q) = 1$, $0<1/q< 1/q_d$, where $q_d = (d^2+d+2)/2$.
Let us fix a number $R>1$ and set $I_R=[0, R]$ and let
$$T_R f(x) = \int_0^R f(t) w(t) e^{-i\inn{x}{\gamma(t)}} dt.$$
It suffices to show
\begin{equation} \label{estimate} \|T_R f\|_{L^q(\bbR^d)}\le C\|f \|_{L^{p}(w dt)}
\end{equation}
with a constant $C$ independent of $R > 1$.   We need to prove this
for
$2<p<q_d$, $q=d(d+1)p'/2$  (and $q_d>2$ if $d\ge 3$); the estimate for $p\le 2$ follows then by interpolation with the trivial case $p=1$.

Observe that \eqref{estimate} holds with some $C=C(a, R)<\infty$, by
the estimates for the nondegenerate curve $\gamma$ (restricted to $I_R$);
notice that indeed $|w(t)|\ge C(a) \min\{1,
e^{R (\sum_j a_j)2/(d^2+d)}\} >0$ on $I_R$.
Let now $\mathcal B_{R,a}$
be the infimum over all $C$ for which \eqref{estimate} holds.
 $\mathcal B_{R,a}$ is finite and we have to establish that
$\mathcal B_{R,a}$ is
uniformly bounded in $R\ge 1$ and $a=(a_1,\dots, a_d)$.

We shall estimate the
$d$-linear  expression
$$\prod_{j=1}^d T_R f_j(x)
= \int_{I_R^d}e^{i \inn{x}{\gamma(t_1)+ \cdots +\gamma(t_d)}}
   \prod_{j=1}^d \left[ f_j(t_j) w(t_j)\right] dt_1 \cdots dt_d.$$
We change variables  $\kappa_{j}(h)=\sum_{i=1}^{j-1}h_i$ as in \eqref{kappah},
 and let
$J_R$ denote the set of all
$h\in [0,R]^{d-1}$ satisfying $\ka_d(h)\le R$. For $h\in J_R$ let
$I_{R,h}= [0, R-\kappa_{d}(h)]$, and  define for any permutation $\pi$ on $d$ letters
\begin{equation}
 \label{Fpi}
F^\pi(h,t)=\chi_{J_R}(h) \chi_{I_{R,h}}(t) \prod_{i=1}^d
f_{\pi(i)}(t+\kappa_i(h)).
\end{equation}
For fixed $h$ let
\begin{equation}\label{Gath}
\Gamma(t,h)=\sum_{j=1}^d \gamma(t+\ka_j(h))
\end{equation}
and
\begin{equation}\label{Hth}H(t,h)=\prod_{j=1}^d w(t+\ka_j(h)).
\end{equation}
Define an operator $S_{R,h}$ by
$$S_{R,h} [F](x)
= \chi_{J_R}(h)\int_{I_{R,h}} e^{i \inn{x}{\Gamma(t,h)}} F(t,h) H(t,h) \,dt. $$
Then
\begin{equation} \label{multilin}
\prod_{j=1}^d T_R f_j(x)= \sum_{\pi\in \fS^d} \int
S_{R,h} [F^\pi](x) dh.
\end{equation}

We first give an estimate on the operators $S_{R,h}$ which will put us in the position to apply the Vandermonde estimate \eqref{vandineq}.

\begin{proposition} \label{DMprop}  Fix $1<p< q_d=\tfrac{d^2+d+2}2$ and let $q=\tfrac{d(d+1)}2 p'$.
For $\vth \in (0,1)$ define
\begin{equation}\label{vth}
\begin{aligned}
&\frac 1{A}=  1-\frac \vartheta {2}, \quad
&&\frac{1}{B}= \frac{1}{p}+\vth(\frac{1}{2}-\frac {1}{p}),
\\
&\frac{1}{s}=
\frac{1-\vartheta}{q}+\frac{\vartheta}{2}, \quad
&&\eta =1- \frac{d+1}{2q}(1-\vth).
\end{aligned}
\end{equation}
Then (with $v$ as in
\eqref{vh})
\begin{multline}\label{SRest}
\Big\|\int S_{R,h}[ F] dh \Big\|_s \\
\le C \cB_{R,a}^{1-\vth}
\Big(\int \Big(\int \big|
F(t,h) H(t,h)^{\eta-\frac{d+1}{4}\vth}\big|^{B}dt
\Big)^{\frac{A}{B}}
 v(h)^{1-A} dh\Big )^{\frac 1{A}}.
\end{multline}
\end{proposition}

\begin{proof} The proof relies on arguments in the papers by
Drury and Marshall  \cite{DM1}, \cite{DM2}. We begin with a few remarks on the
affine arclength measure for the curve $\gamma$ and for the
``offspring''  curves $t\mapsto \Gamma(t,h)$.
Let $\tau$ and $w$ be as in \eqref{taudef}, \eqref{wdef} (for the curve
$\gamma$ in \eqref{expcurv}).
Then
$$            |\tau(t)|=
v(a) \exp\big(t\sum_{i=1}^d a_i\big)$$
with $a=(a_1,\dots, a_d)$,
and
$$ H(t,h)^{1/d} = w(t) \exp\big(\tfrac{2}{d^2(d+1)}
(\sum_{i=1}^d  a_i)
(\sum_{j=2}^d \ka_j(h))
\big) .
$$

Next, $\Gamma(t, h)= \gamma(t) E(h)$, where $E(h)$ is a
$d\times d$ diagonal matrix with the diagonal entries
$$E_{ii}(h)=\sum_{j=1}^d e^{a_i \ka_j(h)}$$
so that $1\le E_{ii}(h)\le d$, for $1\le i\le m$
and $e^{a_i \ka_d(h)} \le E_{ii} (h)\le d e^{a_i \ka_d(h)}$, for
$m+1 \le j\le d$ (for the definition of $m$ see the paragraph after \eqref{expcurv}).
Moreover, if $\tau_h$ is the expression \eqref{taudef} for the curve
$\Gamma(\cdot,h)$ then
$$\tau_{h}(t)= v(a) \exp\big(t\sum_{i=1}^d a_i\big) \prod_{k=1}^d E_{kk}(h).$$

We first   establish the inequality
\begin{equation}\label{uniform} \big\|  S_{R,h} [F H^{-\frac{d-1} d}]\big\|_q
\le C \mathcal B_{R,a} \big( \int |F(t,h)|^p
H(t,h)^{1/d}
dt \big)^{1/p}
\end{equation}
with a constant $C$ uniform in $h$. Notice that the quotient of
$H(t,h)^{1/d} $ and $w_h: =\tau_h^{2/(d^2+d)}$ is independent of $t$, namely
$$Q(h):= \frac{H(t,h)^{1/d}}{w_{h}(t)}=
\frac{
\exp\big(\tfrac{2}{d^2(d+1)}
(\sum_{j=1}^d  a_j)
(\sum_{j=2}^d \ka_j(h))
\big) }
{\big(\prod_{i=1}^d (\sum_{j=1}^d e^{a_i\ka_j(h)})\big)^{\frac{2}{d(d+1)}
}}.
$$
Since $\Gamma(t, h)= \gamma(t) E(h)$ we have by affine invariance
$$\Big(\int \big|\int_{I_{R,h}} e^{i \inn{x}{\Gamma(t,h)}} g(t)w_h(t) dt
\big|^q dx\Big)^{1/q} \le \cB_{R,a} \Big(\int|g(t)|^p w_h(t) dt\Big)^{1/p}
$$
and thus with $g(t):= F(t,h)$,
\[ \big\| S_{R,h}[F H^{-(d-1)/d}]\big\|_q
\le  \mathcal B_{R,a} Q(h)^{1-1/p}
\Big( \int |F(t,h)|^p
H(t,h)^{1/d}dt \Big)^{1/p}.\]
Thus,
the estimate
 \eqref{uniform} will follow once we establish the inequality that
$Q(h)$ is bounded.
But note that
\[\begin{aligned}
Q(h)^{d(d+1)/2} &\le \exp\big(\tfrac{1}{d}
(\sum_{j=1}^m  a_j)
(\sum_{j=2}^d \ka_j(h))
\big)
\frac{
\exp\big(\tfrac{1}{d}
(\sum_{j=m+1}^d  a_j)
(\sum_{j=1}^d \ka_j(h))
\big)}
{\prod_{i= m+1}^d
e^{a_i \ka_d(h)}}
\\
&\le\exp\big(\tfrac{1}{d}
(\sum_{j=m+1}^d  a_j)
(\sum_{j=1}^d (\ka_j(h)-\ka_d(h)))
\big) \le 1
\end{aligned}
\]
since $\ka_d\ge \ka_{d-1} \ge \ka_2\ge \ka_1=0$ and $a_i<0$ for $i\le m$,
$a_i>0$ for $i>m$.
Thus \eqref{uniform} is proved.

We may replace  $F$ by $FH^{(d-1)/d}$ and integrate the resulting estimate
with respect to $h$.
This yields
\begin{equation} \label{pqbound}
 \int \big \| S_{R,h}[F]\big\|_q dh
\le C \mathcal B_{R,a}\int \Big(\int
|F(t,h) H(t,h)^{\tfrac{d-1}d + \tfrac{1}{dp}}|^p dt\Big)^{1/p} dh .
\end{equation}
Note that this implies the claimed estimate \eqref{SRest}
for the case $\vth=0$.

Now as in \cite{DM1}, \cite{DM2} one can
perform  the change of variables $(t,h) \mapsto
\Gamma(t,h)$ (justified in \cite{DM2}, p. 549)
and use Plancherel's theorem, to  obtain
\begin{equation} \label{planchappl}
\Big\|\int S_{R,h}[F]dh\Big\|_2
\le C \Big( \iint \big|F(t,h)H(t,h)J(t,h)^{-1/2}\big|^2 dt
\, dh \Big)^{1/2}
\end{equation}
where $J(t,h)$ is the Jacobian of this transformation.

Interpolating these two estimates gives
\begin{multline}\label{ineq}
\Big\| \int S_{R,h}[F] dh\Big\|_{s}
\\
\le C  \mathcal B_{R,a}^{1-\vth} \Big(\int \Big(\int
\big|F(t,h)H(t,h)^{ \eta } J(t, h )^{-\vartheta/2 }\big|^{B(\vth)}
dt\Big)^{A(\vth)/B(\vth)} dh \Big)^{1/A(\vth)}
\end{multline}
where $0\le \vartheta \le 1$ and $A,B,s,\eta$ are as in \eqref{vth}.
\medskip

We now use a crucial estimate  concerning the determinant of the
$d\times d$  matrix
$\mathcal E(a,s):= \big( e^{a_i s_j}\big)_{i,j=1,\dots d}$.

\medskip

\noi{\bf A total positivity bound by Drury and Marshall.}
\textit{\cite{DM2}, p.546}.
\textit{The estimate}
\begin{equation}\label {tp}
\frac{\det \mathcal E(a,s)}{\prod_{1\le i<j\le d}\big((a_j-a_i)(s_j-s_i)\big)}
\ge c_d \exp\big( \frac 1d \big(\sum_{j=1}^d a_j\big)\big(\sum_{j=1}^d s_j\big)\big)
\end{equation}
{\it holds for all real $a_1,\dots, a_d$ and all real $s_1,\dots, s_d$ with a
constant $c_d$ that depends only on the dimension $d$}.

\medskip

This means $ J(t,h) \ge c_d \ v(h) H(t,h)^{(d+1)/2}$
and therefore
%
\begin{multline*}  \Big\|\int S_{R,h}[F] dh\Big\|_s
 \le C ( \mathcal B_{R,a})^{1-\vartheta } \times \\
\Big(\int \Big(\int
\big| F(t,h) H(t,h)^{\eta- \frac{d+1}{4}\vth}
 \big|^B dt\Big)^{A/B}
v(h)^{-\vartheta A/2}dh \Big)^{1/A}.
\end{multline*}
Now observe  that $A^{-1}=1-\vth/2$ means
$-\vartheta A/2=1-A$ and thus the  proof of the proposition is complete.

\end{proof}

\medskip

\noi {\it Proof of Theorem \ref{powerthm}, continued.}
Proposition \ref{DMprop} enables
us to apply the inequality \eqref{vandineq}.
We wish to use it for the value
\begin{equation}\label{vthp}
 \vth=\vth(p)= \frac{4(d-1)}{(d+1)d p'-4} = \frac{2(d-1)}{q-2}
\end{equation}
and we let $A=A_p$, $B=B_p$, $s=s_p$ and $\eta=\eta_p$ be the values which correspond to $\vth=\vth(p)$ via \eqref{vth}.
The reason for this choice is that the  exponent of $H$ in \eqref{SRest} becomes
\begin{equation}\label{etap}\eta_p -\frac{d+1}{4}\vth= \frac 1p;
\end{equation}
moreover
\begin{equation}\label{sp} s_p = \frac{q}{d}= \frac{d+1}{2} p' .
\end{equation}

In order to apply \eqref{vandineq} we  need
the additional restriction  $1<A_p<\frac{d+2}d$, which corresponds to
$\vth(p)< 4/(d+2)$. A short calculation reveals that this
 requirement is equivalent with our assumption
 $p<\frac{d^2+d+2}2$.

We also set $\sigma_p=2/(d+2-dA_p)$ and obtain after a short computation
\begin{equation*}
\frac{1}{A_p\sigma_p}= \frac{d+2}{A_p}- \frac{d}{2} =
\frac{p^{-1}-q^{-1}}{ 1- 2 q^{-1}}
\end{equation*}
and
\begin{equation*} \frac {1} {B_p} = \frac {1}{p} + \frac{
(\frac 12-\frac 1p)\frac{d-1}{dp'} }
{\frac{d+1}4-\frac{1}{dp'}}.
\end{equation*}
We check that $B^{-1}>p^{-1} >(A_p\sigma_p)^{-1}$  since $2<p$ and we have
\begin{equation} \label{consistency}
 \frac{d-1}{A_p\sigma_p} + \frac {1}{B_p}=\frac dp.
\end{equation}
Now let $\Sigma(A_p, B_p)$ be the simplex defined in the statement of
Proposition
\ref{vandthm}. We apply this proposition
 to  the right hand side of \eqref{SRest}
with $F=F^\pi$ as in \eqref{Fpi}; then by \eqref{multilin}
\begin{equation}
\label{dlinear}
\Big\| \prod_{j=1}^{d} T_R f_j \Big\|_{q/d} \leq C
\mathcal B_{R,a}^{1-\vth(p)} \prod_{j=1}^d \|f_j w^{1/p}\|_{p_j ,1},
\end{equation}
for all $(p_1^{-1},\dots, p_d^{-1})\in \Sigma(A_p,B_p)$.

We continue to argue as in the proof of Theorem \ref{osc} and consider now
the $q_d$-linear expression
$$\fM_R [g_1,\dots, g_{q_d}] =\prod_{k=1}^{q_d} T_R[ g_j w^{-1/p}].$$

Consider the set $K_p$ consisting of the points
$P^\vpi=(P^\vpi_1,\dots, P^\vpi_{q_d})$, $\vpi \in \fS_{q_d}$
({\it i.e.} a permutation on $\{1,\dots, q_d\}$)  with $P^\vpi$
defined by
\[P^\vpi_{\vpi(i)}=
\begin{cases}  1/B_p, &i=1,
\\
1/(A_p\sigma_p),\quad &2\le i\le d,
\\
1/p, &d+1\le i\le q_d.
\end{cases}
\]
The (closed) convex hull of $K_p$
is a simplex on the hyperplane
$\{X\in \bbR^{q_d}:\sum_{i=1}^{q_d} X_i=q_d/p\}$ with vertices
 $P^{\vpi}$, and center $(p^{-1},\dots, p^{-1})$.

By \eqref{dlinear} and H\"older's inequality
\begin{equation*}
\big\|
\fM_R [g_1,\dots, g_{q_d}]\big\|_{q/q_d}
\le C (\cB_{R,a})^{(1-\vth)q_d/d} \prod_{k=1}^{q_d} \|g_k\|_{p_k,1}
\end{equation*}
for all  $(p_1^{-1}, \dots, p_{q_d}^{-1})\in K_p$.
We now apply Proposition \ref{trick} and observe that since $q>q_d$ our
multilinear operator takes values in a Banach space.

Thus we get
$$\big\|
\fM_R [g_1,\dots, g_{q_d}]\big\|_{q/q_d}
\le C (\cB_{R,a})^{(1-\vth)
q_d/d}
 \prod_{k=1}^{q_d}
\|g_{\vpi(k)}\|_{p_k,q_d},$$
for all
$(p_1^{-1}, \dots, p_{q_d}^{-1})\in (\conv K_p)^o$.
Clearly the center
$(p^{-1},\dots, p^{-1})$ belongs to $(\conv K_p)^o$
and it follows that for
$g_i= fw^{1/p}$,
\begin{equation*}
\|T_R f\|_q=
\big\|
\fM_R [g_1,\dots, g_{q_d}]\big\|_{q/q_d}^{1/q_d}\le
C(p,d) (\cB_{R,a})^{(1-\vth)/d}
\|f w^{1/p}\|_{p,q_d}.
\end{equation*}
By $p<q_d$
and the continuous imbedding $L^p\subset L^{p,q_d}$
we have
$$\|f w^{1/p}\|_{L^{p,q_d}} \le
\|f w^{1/p}\|_{L^{p}} =\Big(\int |f(t)|^p w(t) dt\Big)^{1/p}.
$$
Thus $\cB_{R,a}\le C'(p,d) (\cB_{R,a})^{(1-\vth)/d}  $ and the assertion of the theorem follows. \qed

\section{Proof of Theorem \ref{BOext}}\label{sectBOext}
Let $\alpha<\beta$, $\alpha, \beta\notin \{0,1\}$.
We first note that the affine arclength measure $w(t) dt$   for the curve
$(t, t^\alpha, t^\beta)$, $t>0$,
 is given via  $w(t)=c(\alpha,\beta) t^{(\alpha+\beta-5)/6}$
with $c(\alpha,\beta)^6= \alpha\beta(\alpha-1)(\beta-1)(\beta-\alpha)$.
We consider the  case $\alpha+\beta=5$ which clearly plays a
special role as the  affine arclength measure  is now a constant multiple of
Lebesgue measure on $\bbR$.
The case $\alpha=2$, $\beta=3$ has been handled in \S\ref{lorentzbd}, and
part (i) of Theorem   \ref{BOext} asserts that it holds also true for
$\alpha=5-\beta<2$.

To prove this assertion we  consider a more general class of  curves
\begin{equation} \label{yz}
t\mapsto (t, y(t), z(t)), t\in I=(0,b); \quad 0<b < \infty
\end{equation}
where $y,z\in C^3(I)$ and
satisfy a strong nondegeneracy condition introduced in \cite{BO},
namely
\begin{equation} \label{strongnd}
\Delta(s,t):=\big| y''(s)z'''(t)- y'''(s) z''(t)\big|\ge \delta>0, \quad s,t,\in I;
\end{equation}
moreover it is assumed that
\begin{equation}\label {z3}
z'''(t)\neq 0, \quad t\in (0,b),
\end{equation}
however no upper bounds for the third derivatives are  required on the open interval $(0,b)$.
Note that the  determinant in \eqref{strongnd}
cannot change sign. In particular, if
 $h_1, h_2\ge 0$,
$h_1+h_2<b$, and if we consider
the offspring curves
$\Gamma(t,h)=\tfrac 13\sum_{i=1}^3 \gamma(t+\ka_i(h))$, $t<b-h_1-h_2$,
then
$\Gamma(t,h)=(t, y_h(t), z_h(t))$ where
$(y_h,z_h)$ satisfies \eqref{strongnd} (with the same $\delta$)
 on the interval $(0, b-h_1-h_2)$.
This follows from an expansion using the
multilinearity of the determinant.

Let
\begin{equation*}\label{}
\cE f(x):=  \int_0^b e^{-i \inn{x}{ \gamma(t)}} f(t) dt .
\end{equation*}

\begin{proposition} Let $\gamma$ be as in
\eqref{yz}, \eqref{strongnd}.
Then
\begin{equation}\label{Ef}
\| \cE f \|_{L^{7,\infty}(\bbR^3)} \le C \delta^{-1/7} \| f\|_{L^7(I)}
\end{equation}
\end{proposition}

\begin{proof}
Let $\cK(b,\delta)$ be the class of curves $\gamma$ satisfying
\eqref{yz}, \eqref{strongnd} on $(0,b)$ and let
\begin{multline} \label{affrestr4}
\cA_{\delta}(b,R) := \sup_{\rho>0} (1+ R^{-1}\rho)^{-2d}
\sup_{\substack{
\gamma\in \cK(a,\delta)\\0<a \le b
} }
\sup_{\substack{ \|f\|_{L^{7}(0,a)}\le 1}}
\big\|\cE f\big\|_{L^{7,\infty}(B_R)}.
\end{multline}
Clearly $\cA_{\delta}(a,R)
\le C(b,R)<\infty$ for $a\le b$
and we need to show that $\cA_{\delta}(a, R)$ is uniformly bounded in $a$ and $R$.

Now let $a\le b$ and let $\gamma\in \cK(a,\delta)$.
As in \S\ref{osc}  we estimate the trilinear expression
 $\cM(f_1,f_2,f_3)=\prod \cE f_i(x)$ and split $\cM=\sum_k \cM_k$
where
$$\cM_k f(x) = \int_{S_k} e^{-i\inn{x}{ \sum_{i=1}^3\gamma(t_i)}}
\prod f_i(t_i) dt_1dt_2dt_3
$$
with $S_k=\{(t_1,t_2,t_3): 2^{-k-1}<V_3(t)\le 2^{-k}\}$.

It was observed in Lemma 2 of \cite{BO}
that the map $(t_1,t_2,t_3)\to \tfrac 13\sum_{i=1}^3\gamma(t_i)$ is one-to-one, when restricted to $\{t_1<t_2<t_3\}$ (this uses \eqref{z3}).
Denote the  Jacobian of this mapping  by
 $J(t_1,t_2, t_3)$.
Also as in \cite{D}, \cite{BO} one uses a generalized mean value theorem
(\cite{ps}, V.1.95) to obtain the inequality
$J\ge \delta V$. As before this leads to the $L^2$ bound
\begin{equation}\label{EL2}
\| \cM_{k} (f_1, f_2, f_3)\|_2 \le C\delta^{-1/2}
2^{k/6} \| f_1\|_2 \| f_2\|_{\infty} \| f_3\|_{\infty} .
\end{equation}
On the other hand, applying the definition of $\cA$ to the off-spring
curve and the fact that the measure of
$\{(t_2,t_3): (t_1, t_2, t_3)\in S_k\}$ is $O(2^{-2k/3})$ leads to
\begin{equation}\label{E7}
\| \cM_{k} (f_1, f_2, f_3)\|_{L^{7,\infty}(B_R)} \le C
\cA_{\delta}(b,3R)
2^{-2k/3}  \| f_1\|_{7} \| f_2\|_{\infty} \|
f_3\|_{\infty} .
\end{equation}

From here on we argue as in the proof of Theorem \ref{osc}.
Applying Bourgain's interpolation lemma  we get
\begin{equation*}\label{}
\| \cM (f_1, f_2, f_3)\|_{L^{7/3,\infty}(B_R)} \le C \delta^{-2/5}
\cA_{\delta}(b,3R)^{1/5}
  \| f_1\|_{7/3,1} \| f_2\|_{\infty} \|
f_3\|_{\infty}
\end{equation*}
and applying the multilinear interpolation arguments to the symmetric
$n$-linear
expression $\prod_{i=1}^n \cE f_i$, for $n>7$ ({\it e.g.}   $n=21$ as in
\S\ref{lorentzbd}) yields
\begin{equation*}\label{}
\Big\| \prod_{i=1}^n \cE f_i\Big\|_{L^{7/n,\infty}(B_R)} \le C \delta^{-2n/15}
\cA_{\delta}(b, 3R )^{n/15} \prod_{i=1}^n \| f_i\|_{L^{7,r_i}(I)}
\end{equation*}
where $\sum_{i=1}^n r_i^{-1}=n/7$. We may set $f_i=f$, $r_i=7$ and obtain
\begin{equation*}\label{}
\cA_{\delta}(b,R)
\le C \delta^{-2/15}  \cA_{\delta}(b,3R)^{1/15}
\end{equation*}
and since from definition
\eqref{affrestr4} it follows  that
$\cA_{\delta}(b,3R) \lc \cA_{\delta}(b,R)$ we obtain
$\cA_{\delta}(b,R)\lc \delta^{-1/7}$ which is the assertion.
\end{proof}

\begin{proof}[Conclusion of the proof of Theorem \ref{BOext}]
We first consider part  (i). By symmetry we may assume $\alpha\le 2$.
 The cases $\alpha=0$ and $\alpha=1$ are
trivial since then $w_\gamma\equiv 0$, and the case $\alpha=2$
has been already handled in \S\ref{lorentzbd}. Thus suppose
$\alpha=5-\beta<2$, and
$\alpha\neq \{0,1\}$.
Then  by the discussion in the beginning of this section
the affine arclength measure is $c_\alpha^{1/6} (5-2\alpha)^{1/6}dt$
with $c_\alpha=
|\alpha(5-\alpha)(\alpha-1)(4-\alpha)|$.
Moreover
$
\Delta(s,t)=c_\alpha |
(3-\alpha) s^{\alpha-2}t^{2-\alpha}+(2-\alpha)s^{3-\alpha} t^{\alpha-3} |
$
and $c_\alpha^{-1}\Delta(s,t)$ has its minimum $(5-2\alpha)$ at $(1,1)$.
From this part (i) of the theorem follows easily. Part (ii) follows from part (i) by the change of variable $u=t^\alpha$, and interchanging the first and second components of $\gamma$.
\end{proof}

\appendix
\section{\\
Vandermonde operators:
Proof of the Drury-Marshall bound}\label{app2}

For the sake of self-containedness
we give the full proof of Proposition  \ref{vandthm},
due to Drury and Marshall. This is done by first checking
(i) for $d=2$ and $d=3$, and then by  arguing by induction, applying a
special case of  (ii) in $d-2$ dimensions to prove (i) and (ii) in $d$
dimensions.

We note that by a homogeneity argument it suffices to prove that the set
$\Omega_d(1)$ has finite measure in $\bbR^{d-1}$. It is obvious that the measure of $\Omega_2(1)$ is equal to $1$. If $d=3$ then
   $v_3(h)=V_3(\kappa(h_1,h_2))=h_1h_2(h_1+h_2)$
and the set $\{h:v_3(h)\le 1\}$ is contained in the union of two sets
$E_1\cup E_2$   where
$E_1=\{h\in (0,\infty)^2: h_1h_2^2\le 1, h_1\le h_2\}$ and
$E_2=\{(h_1,h_2): (h_2,h_1)\in E_1\}$.
Both sets have area equal to $\int_0^\infty \min\{s, s^{-2}\} ds=3/2$.

Now we  assume that  (ii) has been established in all
dimensions $\le d-1$, and we shall prove that
 $\Omega_d(1)$ has finite measure in $(\bbR_+)^{d-1}$,
and that \eqref{vandineq} holds in $d$ dimensions.

We now set $r\equiv r(h)=\ka_d(h)=h_1+\dots+ h_{d-1}$,
$t\equiv t(h)=h_1/\kappa_d(h)$ and $\tau_i(h)= h_{i+1}/\kappa_d(h)$, $i=1,\dots, d-3$.
We use the change of variable $h\mapsto (t,\tau_1,\dots, \tau_{d-3}, r)$
and observe  the determinant of its derivative is $r(h)^{-d+2}$.

Set $\tilde\ka(\tau)= (\tilde \ka_1,\dots,\tilde \ka_{d-3})$,
with $\tilde \kappa_1(\tau)=0$ and
$\tilde \kappa_i(\tau)=\sum_{k=1}^{i-1} \tau_k$, for $2\le i\le d-2$.
Then we can write
\begin{equation*}\begin{aligned}
v_d&(h)= \Big(\prod_{j=2}^d \kappa_j(h)\Big)
\Big(\prod_{2\le i<j\le d-1} \kappa_j(h)-\kappa_i(h)\Big)
\Big(\prod_{k=2}^{d-1} (\kappa_d(h)-\kappa_{i}(h))\Big)
\\&=\kappa_d(h)^{\frac {d(d-1)}{2}} \Big(\prod_{j=2}^d
\frac{\kappa_j(h)}{\kappa_d(h)}
\Big(1- \frac{\kappa_j(h)}{\kappa_d(h)} \Big)\Big)
\Big(\prod_{2\le i<j\le d-1} \frac{\kappa_j(h)-\kappa_i(h)}{\kappa_d(h)}\Big)
\\&=
r(h)^{\frac{d(d-1)}{2}} \prod_{j=1}^{d-2} [(t+\tilde \kappa_{j}(\tau))
(1-t-\tilde \kappa_j(\tau))]\Big( \prod_{1\le i<j\le d-2}(\tilde \kappa_j(\tau)-\tilde\kappa_i(\tau))
\Big).
\end{aligned}
\end{equation*}
Thus, if $U(s):=\big(\prod_{i=1}^{d-2} [s_i(1-s_i)]\big)V_{d-2}(s)$, defined on  $\Sigma_{d-2}:=\{ s\in \bbR^{d-2}:
0\le s_1\le...\le s_{d-2}\le 1\}$,  then
\begin{equation*}\begin{aligned}
|\Omega_d(1)|&\le \int_{\Sigma_{d-2}}
\int_0^{U(s)^{-\frac{2}{d(d-1)}}}
r^{d-2} dr ds
\le \int_{\Sigma_{d-2}}U(s)^{-2/d} ds
\\
&=
\int_{(\bbR_+)^{d-3}} \int_0^\infty \prod_{i=1}^{d-2}
 |g(t+\tilde\ka_i(\tau))|^A dt\,[v_{d-2}(\tau)]^{1-A} d\tau;
\end{aligned}
\end{equation*}
here $v_{d-2}(\tau):= V_{d-2}(\tilde \kappa(\tau))$,  $A=(d+2)/d$,
$g(s)={s(1-s)}^{-2/(d+2)}\chi_{[0,1]}(s)$.
Thus the last expression is the $A^{\text{th}}$ power of the $L_v^A(L^A)$
 norm of the Vandermonde operator in $d-2$ dimension,
applied to the functions $f_i=g$, $i=1,\dots, d-2$. The value $A=(d+2)/d<
d/(d-2)$ is permissible for the application of part (ii) in $(d-2)$ dimensions.
Now with $\sigma=2(d-(d-2)A)^{-1}=d/2$ we need to verify that
$g\in L^{p,1}$ with $(d-2)/p=A^{-1}(1+(d-3)/\sigma)$
(which corresponds to the point in the center of the $(d-3)$-dimensional
simplex $\Sigma(A,A)$).
Note that $p=(d+2)/3$, and as $g$ belongs to $L^r[0,1]$ for all $r<(d+2)/2$
it belongs surely to $L^{p,1}$. Thus part (i) is verified in $d$ dimensions.

We now turn to the proof of \eqref{vandineq} in $d$ dimensions.
First notice that the allowable $p_i$'s are given by the equation
$\sum_{i=1}^d p_i^{-1} = (d-1)(\sigma A)^{-1}+B^{-1}$
and that $(p_1^{-1}, \dots, p_d^{-1})$ belongs to  $\Sigma(A,B)$ if and
only if
$$\frac 1{p_i} =
\frac {1}{\sigma A}+ \big(\frac 1B -\frac 1{\sigma A}\big) \frac 1{r_i}$$
where $r_i\in [1,\infty]$ with $\sum_{i=1}^d r_i^{-1}=1$.

It suffices to prove the estimate \eqref{vandineq}
 for  $f_i$ which is pointwise dominated
 by characteristic functions of measurable sets $E_i$, $i=1,\dots, d$,
and by monotonicity properties of the operator we  may assume that
$f_i=\chi_{E_i}$. Thus we need to prove
\begin{equation}\label{vandineqchar}
\big\|\frak V(\chi_{E_1},\dots, \chi_{E_d})\big\|_{L^A_v(L^B)}\le
C\prod_{i=1}^d |E_i|^{1/p_i}.
\end{equation}
We use a duality argument for the $h$ integral in
\eqref{mixednormforV}. By part (i)  the function $|v|^{1-A}$
belongs to $L^{\sigma ',\infty}(\Bbb R^{d-1})$
for $\sigma'=2d^{-1}(A-1)^{-1}$. Note that because of our assumption on
$A$  we have $\sigma'\in (1,\infty)$; moreover $\sigma'$ is the  conjugate exponent  to
$\sigma =2/(2+d-dA)$.
Define
$$\Phi(h)= \int \prod_{i=1}^d |\chi_{E_i}(t+\kappa_i(h))| dt,$$
as a function defined on $(0,\infty)^{d-1}$.
As $\chi_{E_i}$ assumes only the values  $0$ and $1$ it suffices to show
that
$\big \|\Phi ^{A/B}\big\|_{L^{\sigma ,1} }
\lc\prod_{i=1}^d |E_i|^{A/p_i}$
which follows from
\begin{equation} \label{phibd}
\|\Phi\|_{L^{A\sigma /B,A/B} }
\le C\prod_{i=1}^d
|E_i|^{\frac {B}{p_i}}.
\end{equation}
We now use the familiar inequality
\begin{equation}\label{convexity}
\|G\|_{L^{P,s}} \le C(P,s) \|G\|_{1}^{1/P} \|G\|_\infty^{1-1/P}
\end{equation}
which holds for all $P\in (1,\infty)$ and all $s\in (0,\infty)$;
and we apply this for $G\equiv \Phi$, and $P= A\sigma /B>1$,
$s= A/B\in (0,1]$.
It is easy to see that
$
\|\Phi \|_1 \le \prod_{i=1}^d|E_i|
$
and, by H\"older's inequality,
$
\|\Phi \|_\infty \le \prod_{i=1}^d|E_i|^{1/r_i}
$
for all $r_1,\dots, r_d\in [1,\infty]$
satisfying $\sum_{i=1}^d r_i^{-1}=1$.
Now
by \eqref{convexity},
\begin{equation} \label{ribd}
\big \|\Phi\big\|_{L^{A\sigma /B,A/B} } \le c(A,B)
\prod_{i=1}^d|E_i|^{ \frac{B}{A\sigma }+(1-\frac B{A\sigma })
\frac 1{r_i}}
\end{equation}
where $r_i\in [1,\infty]$ with $\sum_{i=1}^d r_i^{-1} =1$.
By the above description of the simplex $\Sigma(A,B)$ this inequality
yields \eqref{phibd} and thus the assertion.\qed



\end{document}